\documentclass[11pt]{amsart}
\usepackage{geometry}                
\usepackage{graphicx,placeins}
\usepackage{amssymb}
\usepackage{mathtools}
\mathtoolsset{showonlyrefs}

\def\uv{{\bf u}}

\def\epsit{{{\bf{\epsilon}}}}

\newcommand{\dens}{\rho}
\newcommand{\Hd}{(H^1(\Omega))^2}
\newcommand{\ee}{{\bf{\epsilon}}}
\newcommand{\Hgammad}{\mathbf{{\mathcal V}}}

\newcommand{\uu}{{\mathbf{u}}}
\newcommand{\uud}{\mathbf{u_d}}
\newcommand{\cU}{{\mathbb Q}}
\newcommand{\Tau}{\mathcal{T}}

\def\wn2{w_n^2}
\def\ws2{w_s^2}

\def\tn2{t_n^2}
\def\ts2{t_s^2}

\def\vv{\underline{v}}

\def\vv{{\bf v}}
\def\qq{{\bf q}}
\def\Pepsi{\Pi^\epsilon_E}

\usepackage{changes}
\usepackage{todonotes}

\renewcommand{\emph}[1]{{\it #1}}

\usepackage{setspace}

\newtheorem{assumption}{Assumption}[section]
\newtheorem{remark}{Remark}[section]

\title[VEM and Topology Optimization
]{On the Virtual Element Method for Topology Optimization on polygonal meshes: a numerical study}
\thanks{Paola F. Antonietti has been partially supported by the SIR Starting Grant n. RBSI14VT0S ``PolyPDEs: Non-conforming polyhedral finite element methods for the approximation of partial differential equations" funded by  MIUR - Italian Ministry of Education, Universities and Research. 
M. Verani has been partially supported by the Italian research grant  {\sl Prin 2012}  2012HBLYE4  ``Metodologie innovative nella modellistica differenziale numerica'' and by INdAM-GNCS}
\author{P.~F.~Antonietti}
\address{P.~F.~Antonietti: MOX, Dipartimento di Matematica, Politecnico di Milano, Italy}
\email{paola.antonietti@polimi.it}

\author{M.~Bruggi}
\address{M.~Bruggi: Dipartimento di Ingegneria Civile e Ambientale, Politecnico di Milano, Italy}
\email{matteo.bruggi@polimi.it}

\author{S.~Scacchi}
\address{S.~Scacchi: Dipartimento di Matematica, Universit\`a degli Studi di Milano, Italy}
\email{simone.scacchi@unimi.it}

\author{M.~Verani}
\address{M.~Verani: MOX, Dipartimento di Matematica, Politecnico di Milano, Italy}
\email{marco.verani@polimi.it}

\begin{document}
\maketitle

\begin{abstract}
It is well known that the solution of topology optimization problems may be affected both by the geometric properties of the computational mesh, which can steer the minimization process towards local (and non-physical) minima, and by the accuracy of the method employed to discretize the underlying differential problem, which may not be able to correctly capture the physics of the problem. In light of the above remarks, in this paper we consider polygonal meshes and employ the virtual element method (VEM) to solve two classes of paradigmatic  topology optimization problems, one governed by nearly-incompressible and compressible linear elasticity and the other by Stokes equations. Several numerical results show the virtues of our polygonal VEM based approach with respect to more standard methods.
\end{abstract}

\bigskip


\bigskip

\section{Introduction}
The study of numerical methods for the approximation of partial differential equations on polygonal and polyhedral meshes is drawing the attention of an increasing number of researchers (see, e.g., the special issues  {\cite{BeiraoErn_2016,Bellomo-Brezzi-Manzini:2014}} for a recent overview of the different methodologies). 
Among the different proposed methodologies, here we focus on the Virtual Element Method (VEM) which has been introduced in the pioneering paper \cite{VEM-basic:2013} and can be seen as an evolution of the Mimetic Finite Difference method, see, e.g.,~\cite{MFD:book, Lipnikov-Manzini-Shashkov:2014} for a detailed description. Recently, VEM has been analyzed for general elliptic problems \cite{VEM-elliptic:2016, Berrone-SUPG:2016}, linear and nonlinear elasticity \cite{Vem-elasticity:2013,Paulino-et-al:2014,VEM:inelastic:2015}, plate bending \cite{Brezzi-Marini-plate:2013,Chinosi-Marini:2016}, Cahn-Hilliard \cite{Antonietti-Beirao-Scacchi-Verani:2016}, Stokes \cite{Antonietti-Beirao-Mora-Verani:2014, Beirao-Lovadina-Vacca:2015}, Helmholtz \cite{Perugia-Pietra-Russo:2016}, parabolic \cite{Vacca-Beirao:2015}, Steklov eigenvalue \cite{VEM-Steklov:2015}, elliptic eigenvalue \cite{Gardini-Vacca:2016} and discrete fracture networks \cite{VEM-DFN:2014}. In parallel, several different  variants of the VEM have been proposed and analysed: mixed \cite{Brezzi-Falk-Marini:2014,mixed-vem:2016}, discontinuous \cite{VEM-discontinuous}, $H(\text{div})$ and $H(\bold{curl})$-conforming \cite{Hdiv-vem:2016}, hp \cite{Beirao-Chernov-Mascotto-Russo:2016}, serendipity \cite{Serendipity:2016} and nonconforming \cite{Ayuso-Lipnikov-Manzini:2016, Cangiani-Manzini-Sutton:2015, Cangiani-Gyrya-Manzini:2016, vem-cinesi:2016, ncvem-Antonietti-Manzini-Verani:2016}  VEM. 

Such a flourishing research activity  founds an important motivation in the great flexibility that the use of polytopal meshes can ensure in dealing with problems posed on very complicated and possibly deformable geometries. In this respect, as first recognized by G.H. Paulino and his collaborators in a series of ground breaking papers \cite{Talischi-Paulino-Pereira-Menezes:2010, PolyTop, PolyMesher, Paulino-fluid, Gain-Paulino-Duarte-Menezes:2015},  topology optimization represents an intriguing challenge for the use of polyhedral meshes. Topology optimization is a fertile area of research that is mainly concerned with the automatic generation of optimal layouts to solve design problems in Engineering. The classical formulation addresses the problem of finding the best distribution of an isotropic material that minimizes the work of the external loads at equilibrium,
while respecting a constraint on the assigned amount of volume. This is the so-called minimum compliance formulation that can be conveniently employed to achieve stiff truss-like layout within a two-dimensional domain. A classical implementation resorts to the adoption of four node displacement-based finite elements that are coupled with an elementwise discretization of the (unknown) density field. When regular meshes made of square elements are used, well-known numerical instabilities arise, see in particular the so-called checkerboard patterns. On the other hand, when unstructured meshes are needed to cope with complex geometries, additional instabilities can steer the optimizer towards local minima instead of the expected global one. Unstructured meshes approximate the strain energy of truss-like members with an accuracy that is strictly related to the geometrical features of the discretization, thus remarkably affecting the achieved layouts.
On this latter issue, as pointed out also in \cite{Talischi-Paulino-Pereira-Menezes:2010}, the use of polyhedral meshes provide flexibility in the difficult computational task of meshing complex geometries, while, in parallel, it can contribute to avoid that the  geometry of the mesh 
dictates the possible layout of material and the orientation of members, thus excluding physical optimal configurations from the final design obtained by the numerical procedure. 
  
The aim of this paper is to push forward the study of \cite{Talischi-Paulino-Pereira-Menezes:2010, Gain-Paulino-Duarte-Menezes:2015, Paulino-fluid}. In \cite{Talischi-Paulino-Pereira-Menezes:2010} the authors analyze the possibility of avoiding sub-optimal (non-physical) layout in topology optimization for structural applications when {\em polygonal finite elements} and polytopal meshes are employed, whereas in \cite{Gain-Paulino-Duarte-Menezes:2015} the virtual element method is employed for solving compliance minimization and compliant mechanism problems in three dimensions. In \cite{Paulino-fluid} {\em polygonal finite elements} are employed to solve topology optimization problems governed by Stokes equations on polygonal meshes. In view of the above contributions, and with the goal of deepening the  comprehension of the role of VEM and polygonal meshes in topology optimization, we focus on the use of this latter method for solving topology optimization governed by linear elasticity (compressible and nearly-incompressible) and Stokes flow. For each of the above examples, we systematically consider the impact that the combined approach VEM and polygonal meshes has on the quality of the obtained layout and compare them with the ones provided by standard approaches.

The outline of the paper is the following. In Section \ref{S:2} we present the continuous formulation of the topology optimization problems that we will consider throughout the paper, while in Section \ref{S:3} we introduce the corresponding virtual element discretizations. In Section \ref{S:4} we present and extensively discuss several numerical experiments assessing the virtues of the combined use of VEM and polygonal meshes in solving each of the previously introduced topology optimization problems. Finally, in Section \ref{S:5} we draw some conclusion.

\section{Topology Optimization problems: continuous formulation}\label{S:2} 
\noindent We briefly recall the continuous formulations of the topology optimization problems we are interested in, namely the minimum compliance problem governed by the linear elasticity equation (Section~\ref{S:compliance}) and the optimal flow problem governed by the Stokes equation (Section~\ref{S:stokes}). We first recall some notation that will be useful in the following. Let $\Omega$ be a two-dimensional bounded, polygonal domain with boundary $\Gamma=\partial\Omega$ and let $\Gamma_d\subset \Gamma$ be a subset of the boundary of the domain. We introduce the following spaces
\begin{eqnarray}
\Hgammad_0&=&\{\uu\in\Hd:\quad \uu={\bf 0}~\text{on}~\Gamma_d\}\nonumber\\
\Hgammad_d&=&\{\uu\in\Hd:\quad \uu=\uud~\text{on}~\Gamma_d\}\ \nonumber
\end{eqnarray}
where $\uud$ is a possibly null given function.  Moreover, let us introduce the control space 
\[
\mathcal{Q}_{\textrm{ad}}=\{\dens\in L^\infty(\Omega):\quad
0<\rho_{\textrm{min}}\leq\dens\leq 1~\text{a.e.~in~}\Omega\}
\]
of bounded functions representing the material density
in $\Omega$, where $\rho_{\textrm{min}}$ is some positive lower bound.

\subsection{Minimum compliance}\label{S:compliance}
In this section we shortly describe the topology optimization problem for minimum compliance. This corresponds to find the optimal distribution of a given amount of linear elastic isotropic material (described by an element of $\mathcal{Q}_{\textrm{ad}}$) such that the work of the external load against the corresponding displacement at equilibrium is minimized.

More precisely, let $\lambda_0$ and $\mu_0$ be the Lam\'e coefficients of the given material and introduce the bilinear form $a(\cdot,\cdot): \Hd\times \Hd\to \mathbb{R}$ defined as follows
\begin{equation}
a(\uu,{\vv})=2 \mu_0 \int_{\Omega} \ee(\uu):\ee({\vv})~d{\bf x} + 
\lambda_0 \int_{\Omega} \text{div}~\uu~\text{div}~ \vv ~d {\bf x},
\end{equation}
where $\epsit(\uv) =
\frac{1}{2}\left(\nabla\uv+\nabla^T\uv\right)$ is the symmetric gradient.
Moreover, let us introduce the semi-linear form $a(\dens;\cdot,\cdot):\mathcal{Q}_{\textrm{ad}}\times \Hd\times\Hd\to\mathbb{R}$
\begin{eqnarray}
  && a(\dens;\uv,\vv)=2  \int_{\Omega} \mu(\rho) \ee(\uu):\ee(\vv)~d{\bf x} + 
 \int_{\Omega} \lambda (\rho) \text{div}~\uu~\text{div}~ \vv ~d {\bf x}
\end{eqnarray}
where, according to the classical SIMP approach, we set $\lambda(\rho)=\lambda_0 \rho(x)^{p_\lambda}$ and  $\mu(\rho)=\mu_0 \rho(x)^{p_\mu}$, being $p_\lambda$ and $p_\mu$ positive parameters (typically equal to $3$). 
Clearly, when $\rho=1$ in $\Omega$ we have  $ a(\dens;\uv,{\vv})= a(\uv,\vv)$. Finally, let the linear functional $\mathcal{F}(\cdot):\Hd\to\mathbb{R}$ be defined as 
\begin{eqnarray}
  &&\mathcal{F}(\vv)=\int_{\Gamma_t}{\bf f}_{t}\cdot\vv~d{\bf x}\nonumber
\end{eqnarray}
where $\Gamma_t=\Gamma\setminus \Gamma_d$ and the given function ${\bf f}_{t}\in (L^2(\Gamma_t))^2$ represents the external load.
In view of the above definitions, given a material distribution described by the function $\rho$, the elasticity problem (or direct problem) reads as follows: find $\uu \in  \Hgammad_d$ such that  
\begin{equation}\label{pb:elast}
 a(\dens;\uv,\vv)= \mathcal{F}(\vv)
\end{equation}
for any $\vv \in   \Hgammad_0$.
According to the Clapeyron theorem, the continuous formulation of
the topology optimization problem for minimum compliance governed by the elasticity equation can be
therefore written as:
\begin{equation}
  \left\{
  \begin{array}{ll}
    \displaystyle \min_{\dens\in \mathcal{Q}_{\textrm{ad}}} & \mathcal{C}(\dens,\uu) =
    \displaystyle \int_{\Gamma_t}{\bf f}_{t}\cdot\uv~dx \\
    \\
    \mbox{s.t.} & \displaystyle  a(\dens;\uv,\vv)= \mathcal{F}(\vv)
    \mbox{       }\mbox{       }\forall\vv\in\Hgammad_0 \\
    \\
    & \displaystyle\frac{1}{V}\int_{\Omega}\dens dx~ \leq {V_f},\\
  \end{array}
  \right. \label{eq:opti2}
\end{equation}
being ${V_f}$ the available amount of material as a fraction of the
whole domain $V=\int_{\Omega} 1 dx$.  Minimizing the compliance $\mathcal{C}$ of a structure acted upon by a prescribed set of assigned forces means minimizing the work of external loads, i.e. looking for a stiff structure.
 When $\lambda_0\to +\infty$ we refer to  \eqref{eq:opti2} as the minimum compliance problem in the case of nearly-incompressible elasticity, otherwise for moderate values of $\lambda_0$ we refer to it as the minimum compliance problem in the case of compressible elasticity.

\subsection{Optimal Stokes flow}\label{S:stokes} In this section we recall the classical topology optimization problem for optimal Stokes flows  \cite{Borrvall-Petersson:03}: given a design domain $\Omega$ with certain boundary conditions, 
 we are interested in determining at what places of $\Omega$ there should be fluid or 
solid in order to minimize a certain energy functional $\mathcal{E}$ (subject to the constraint of the availability of a given amount of fluid). More precisely, let us introduce the bilinear form
\begin{equation}\label{stokes:form}
a(\dens;\mathbf{u},\mathbf{v})=  2\mu_0 \int_{\Omega} \ee(\mathbf{u}):{\ee}(\mathbf{v})~dx +
\lambda_0 \int_\Omega \text{div}\uv ~\text{div}\vv ~dx + \int_\Omega \alpha(\dens) \mathbf{u} \mathbf{v}
\end{equation} 
where $\mu_0$ is the viscosity of the fluid, $\lambda_0$ is a penalty parameter employed to enforce the incompressibility condition and  $\alpha(\dens)=\frac{5\mu}{\rho^2}$. Given a material distribution described by the function $\rho\in  \mathcal{Q}_{\textrm{ad}} $, the Stokes problem (or direct problem) reads as follows: find $\uu=\uu(\dens) \in   \Hgammad_d$ such that 
\begin{equation}\label{pb_stokes:cont}
a(\dens;\uv,\vv)= 0 \qquad \forall {\vv}\in   \Hgammad_0.
\end{equation}
Thanks to the choice of $\alpha(\rho)$, it turns out that the solution $\uu$ of the above problem is null where 
$\dens=\dens_{\textrm{min}}$ and solves a plane flow model where $\dens=1$ (see  \cite{Borrvall-Petersson:03} for more details). This amounts to interpret the region where $\rho=\rho_{\textrm{min}}$ as
occupied by a solid material, whereas the region where $\dens=1$ as occupied by the fluid.\\

In order to formulate our optimization problem, let us introduce the energy functional $\mathcal{E}(\dens,\mathbf{u}) =
    \displaystyle \frac 1 2 a(\dens;\mathbf{u},\mathbf{u})$ which is a measure of the dissipated energy associated the pair $(\rho,\uu)$. Thus, the optimal flow problem governed by the Stokes equation reads as
\begin{equation}
  \left\{
  \begin{array}{ll}
    \displaystyle \min_{\dens\in \mathcal{Q}_{\textrm{ad}}} & \mathcal{E}(\dens,\mathbf{u})  
    \\
    \mbox{s.t.} & \displaystyle  a(\dens;\uv,\vv)= 0 \qquad \forall {\vv}\in   \Hgammad_0 
    \mbox{       }\mbox{       }\\
    \\
    & \displaystyle\frac{1}{V}\int_{\Omega}\dens dx~ \leq {V_f}.\\
  \end{array}
  \right. \label{eq:opti3}
\end{equation}
\section{VEM discretization}\label{S:3}
In this section we introduce the virtual element discretization of the topology optimization problems 
\eqref{eq:opti2} and \eqref{eq:opti3}, which will be addressed in Sections \ref{S:compliance_d} and 
  \ref{S:Stokes_d}, respectively.  Let $\Tau_h$ represent a decomposition of $\Omega$ into general,  possibly non-convex, polygonal elements $E$
with $\text{diam}(E)=h_E$, where $\text{diam}(E)=\max_{x,y\in E} \| x -y\|$.
In the following, we will denote by $e$ the straight edges of the mesh $\Tau_h$ and, for all $e \in \partial E$, ${\bf n}_E^e$ will denote the unit normal vector to $e$ pointing outward to $E$. We will use the symbol $\mathbb{P}_k(\omega)$ to denote the space of polynomials of degree less than or equal to $k\geq1$ 
living on the set $\omega \subseteq {\mathbb R}^2$. Moreover, we will work under the following mesh regularity  assumption on $\Tau_h$ (see, e.g., \cite{VEM-basic:2013}):
\begin{assumption}\label{meshass}
We assume that there exist positive constants $c_s$ and $c_s'$ such that every element $E\in \{ \Omega_h \}_h$ is star shaped with respect to a ball with radius $\rho \ge c_s h_E$ and every edge $e \in \partial E$ has at least length $h_e \ge c_s' h_E$.
\end{assumption}
Let us first introduce the discrete counterpart of the space $\mathcal{Q}_{\textrm{ad}}$, namely the finite dimensional
space of piecewise constant admissible controls
$$
\cU_{\textrm{ad}}=\{\dens_h\in\mathcal{Q}_{\textrm{ad}}:~{\dens_h}_{\vert
E}\in\mathbb{P}_0(E)~~\forall E \in \Tau_h\}\ .
$$
Clearly, a function $\dens_h\in \cU_{\textrm{ad}}$ is uniquely determined by its value $\rho_E$ in each polygon $E\in \Tau_h$. Hence, the dimension of $\cU_{\textrm{ad}}$ equals the cardinality of $\Tau_h$.

\subsection{Minimum compliance}\label{S:compliance_d}
Following \cite{Vem-elasticity:2013}, it is possible to introduce  the low-order discrete VEM spaces
$\mathbf{V}_{0,h}\subset \mathcal{V}_0$ and $\mathbf{V}_{h}\subset \mathcal{V}$, a discrete form $a_h(\dens_h;{\uu}_h,{\vv}_h)$ approximating $a(\dens;{\uu},{\vv})$   
 and a discrete functional $\mathcal{F}_h({\vv}_h)$ approximating $\mathcal{F}({\vv})$
such that the VEM discretization of \eqref{pb:elast} reads as: given $\dens_h\in \cU_{\textrm{ad}}$ find $\uu_h\in \mathcal{V}_h$ such that
\begin{equation}\label{VEM:pb}
a_h(\dens_h;\uu_h,\vv_h)=\mathcal{F}_h(\vv_h)
\end{equation}
for any $\vv_h\in \mathcal{V}_{0,h}$.
In particular, the global VEM spaces $\mathbf{V}_{0,h}$ and $\mathbf{V}_{h}$ are obtained by gluing suitable local discrete VEM spaces, denoted by $\mathcal{V}_h(E)$, whose elements are uniquely identified by the values at the vertices of the polygon $E$ and contain linear polynomials, i.e $(\mathbb{P}_1(E))^2\subset \mathcal{V}_h(E)$. It is immediate to verify that the dimension of $\mathbf{V}_{0,h}$ (the same happens for $\mathbf{V}_{h}$) equals two times the number of the interior vertices of the partition $\Tau_h$ plus those belonging to $\Gamma_t$, having fixed  the values at vertices belonging to $\Gamma_d$ to incorporate Dirichlet boundary conditions. More precisely, according to \cite{Vem-elasticity:2013}, we have 
$$ 
\mathcal{V}_h(E) =\{{\vv}_h\in (H^1(E))^2: A_{\lambda_0,\mu_0}  {\vv}_h = {\bf 0} \text{~on~} E,   {\vv}_h|_e \in (\mathbb{P}_1(e))^2~~\forall e\in \partial E\}
$$
where $$A_{\lambda_0,\mu_0}  {\uu} = - 
\left( \begin{array}{c}
2\mu_0(u_{1,xx}+\frac 1 2(u_{1,yy}+u_{2,xy}))+\lambda_0(u_{1,xx}+u_{2,yx})\\ 
 2\mu_0(\frac 1 2(u_{1,yx}+u_{2,xx})+u_{2,yy})+\lambda_0(u_{1,xy}+u_{2,yy})
\end{array}
\right).
$$

Let us now define the projection operator 
$\Pepsi:\mathcal{V}_h(E)\rightarrow (\mathbb{P}_1(E))^2$ solution of
\[
\left\{
\begin{array}{ll}
a^E(\Pepsi\vv_h,\qq)=a^E(\vv_h,\qq) & \forall\qq\in (\mathbb{P}_1(E))^2\\
\overline{\Pepsi \vv_h}=\overline{\vv_h} &
\end{array}
\right.
\]
for all $\vv_h\in\mathcal{V}_h(E)$, where $a^E(\cdot,\cdot)$ is the bilinear form $a(\cdot,\cdot)$ restricted to the element $E$
and, for any regular function $\phi$, we set
\[
\overline{\phi}=\frac{1}{n}\sum_{i=1}^n \phi(V_i),\qquad V_i=\mbox{vertices of }E.
\]
It is easy to see that $\Pepsi$ is computable from the degrees of freedom of the local VEM space.\\

The construction of the global form $a_h(\dens_h;\uu_h,\vv_h)$ hinges upon  the construction of local forms $a^E_h(\uu_h,{\vv}_h): \mathcal{V}_h(E) \times \mathcal{V}_h(E)\to \mathbb{R}$  defined as 
\begin{multline*}
a_h^E(\uu_h,{\vv_h})=2 \mu_0 \int_{E}\ee(\Pepsi\uu_h):\ee({\Pepsi\vv}_h)~d{\bf x}  
                     +~\lambda_0 \int_{E}  \text{div}~\Pepsi\uu_h~\text{div}~ {\Pepsi\vv}_h ~d {\bf x}\\
                     +~S^{E,\epsilon}(\uu_h-\Pepsi\uu_h,\vv_h-\Pepsi\vv_h),
\end{multline*}
where the bilinear form $S^{E,\epsilon}$ is a suitable stabilization term with the same scaling properties of the sum
of the first and second term (see  \cite{Vem-elasticity:2013} for more details).
Then, the global form reads as follows:
\begin{equation}\label{form:VEM}
a_h(\dens_h;{\uu}_h,{\vv}_h)=\sum_{E \in {\Tau}_h} \rho_E^p a^E_h(\uu_h,{\vv}_h),
\end{equation}
where we employed the fact that $\rho_h{\vert_ E}\in \mathbb{P}_0(E)$. \\

In the case of nearly-incompressible materials, i.e. for $\lambda_0\to +\infty$, the local discrete bilinear form is modified as follows:
\begin{multline}\label{form_disc:elast}
a^E_h(\uu_h,{\vv}_h)=2\mu_0 \int_{E} \ee({\Pepsi\uu}_h):{\ee}({\Pepsi\vv}_h)~dx + 
                      \lambda_0 \int_E  (\Pi_0 \text{div}\uu_h)(\Pi_0 \text{div}{\vv}_h) ~dx\\
                      +S^{E,\epsilon}(\uu_h-\Pepsi\uu_h,\vv_h-\Pepsi\vv_h),
\end{multline}
where $\Pi_0$ denotes the $L^2$-projection on constant functions and the bilinear form
$S^{E,\epsilon}$ is a suitable stabilization term with the same scaling properties of the first term. 
We refer to  \cite{Vem-elasticity:2013} for more details.
\begin{remark}\label{rem:inf-sup}
For nearly-incompressible elasticity, in \cite{Vem-elasticity:2013} it is not theoretically proved that the resulting lowest-order discrete VEM problem stemming from \eqref{form_disc:elast} is well posed. However, in 
\cite{BeiraoLipni:2010}, it is reported numerical evidence that, 
in the context of Mimetic Finite Differences (MFD), the formulation is well posed on hexagons.
\end{remark}

According to the previous considerations, the VEM discretization of the topology optimization problem \eqref{eq:opti2} reads as follows: 
\begin{equation}
  \left\{
  \begin{array}{ll}
    \displaystyle \min_{\dens_h\in \cU_{\textrm{ad}}} & \mathcal{C}(\dens_h,\mathbf{u}_h) ~=
    \displaystyle   \mathcal{F}_h({\mathbf u}_h) 
    \\
    \mbox{s.t.} & \displaystyle  a_h(\dens_h;{\mathbf u}_h,{\mathbf v}_h)= \mathcal{F}_h(\vv_h)
    \mbox{       }\mbox{       }\forall\vv_h\in\Hgammad_{0,h} \\
    \\
    & \displaystyle\frac{1}{V}\int_{\Omega}\dens_h dx~ \leq {V_f}.\\
  \end{array}
  \right. \label{eq:opti2_h}
\end{equation}

\subsection{Optimal Stokes flow}\label{S:Stokes_d} Hinging upon the results of the previous sections, the virtual discretization of \eqref{eq:opti3} easily follows. Indeed, bearing in mind \eqref{form_disc:elast}, the discrete virtual counterpart of \eqref{pb_stokes:cont} reads as: given $\dens_h\in \cU_{\textrm{ad}}$ find $\uu_h\in \mathcal{V}_h$ such that
\begin{equation}\label{VEM_stoke:pb}
a_h(\dens_h;\uu_h,\vv_h)=\bf{0}
\end{equation}
for any $\vv_h\in \mathcal{V}_{0,h}$, where as usual the global discrete form $a_h(\dens_h;\uu_h,\vv_h)$ is defined in terms of the local forms as  
$a_h(\dens_h;{\uu}_h,{\vv}_h)=\sum_{E \in {\Tau}_h} a^E_h(\uu_h,{\vv}_h) + b^E_h(\dens_h;\uu_h,{\vv}_h)$
with $a^E_h(\uu_h,{\vv}_h)$ being the same as in (\ref{form_disc:elast})
and
\begin{equation}
b^E_h(\dens_h;\uu_h,{\vv}_h)=\int_E \frac{5\mu_0}{\dens_E^2} \Pi^0_E\uu_h \cdot \Pi^0_E\vv_h dx +
S^{E,0}(\uu_h-\Pi^0_E\uu_h,\vv_h-\Pi^0_E\vv_h), 
\end{equation}
where $\Pi^0_E:\mathcal{V}_h(E)\rightarrow (\mathbb{P}_1(E))^2$ is the $L^2$-projection
and the bilinear form $S^{E,0}$, as above, is a suitable stabilization term with the same scaling properties of the
first term.
Note that, using the {\em augmented space argument} (see \cite{Enhanced-VEM:2013}), on the local
VEM space the projections $\Pi^0_E$ and $\Pepsi$ coincide, thus $\Pi^0_E$ is computable.\\

In order to formulate the discrete optimization problem, let us introduce the energy functional $\mathcal{E}_h(\dens_h,\mathbf{u}_h) =
    \displaystyle \frac 1 2 a_h(\dens;\mathbf{u}_h,\mathbf{u}_h)$. Thus, the virtual discretization of \eqref{eq:opti3} reads as
    \begin{equation}
  \left\{
  \begin{array}{ll}
    \displaystyle \min_{\dens\in \cU_{\textrm{ad}}} & \mathcal{E}(\dens_h,\mathbf{u}_h)  
    \\
    \mbox{s.t.} & \displaystyle  a(\dens_h;\uv_h,{\vv}_h)= 0 \qquad \forall {\vv}_h\in   \Hgammad_{0,h} 
    \mbox{       }\mbox{       }\\
    \\
    & \displaystyle\frac{1}{V}\int_{\Omega}\dens_h dx~ \leq {V_f}.\\
  \end{array}
  \right. \label{eq:opti3_h}
\end{equation}

\section{Numerical results}\label{S:4}
Several numerical examples are presented in this section dealing with the VEM discretization of the topology optimization problems introduced above. The minimum compliance problem governed by compressible (Section~\ref{num:compr}) and nearly-incompressible linear elasticity (Section~\ref{num:incompr}) is solved, as well as the optimal flow problem governed by the Stokes equation (Section~\ref{num:stokes}).

The Method of Moving Asymptotes (MMA) \cite{mma}, an algorithm based on sequential convex programming, is herein adopted to tackle the discrete optimization problems \eqref{eq:opti2_h}-\eqref{eq:opti3_h}. 

As described in the previous section, an element--wise density discretization is implemented to approximate the unknown density field. This conventional discrete scheme is affected by well--known numerical instabilities, such as the arising of checkerboard patterns and mesh dependence, see e.g. \cite{4}. Several methods are available in the literature to overcome both problems \cite{sigpet}. Following \cite{torto} and some robust application in stress--constrained topology optimization, see e.g. \cite{filt, 26a,26b,26c}, a filter is herein applied to the density unknowns $\rho_E$ and a new set of physical variables $\tilde{\rho}_E$ is defined as:
\begin{equation}
\displaystyle \tilde{\rho}_E=\frac{1}{\sum_{E'\in\Tau_h} H_{E,E'}}\sum_{E'\in\Tau_h}
H_{E,E'}\rho_{E'},\qquad \quad H_{E,E'}=\sum_{E'\in\Tau_h} {\rm max}(0,r_{\textrm{min}}-{\rm dist}(E,E')).
\label{eq:filter}
\end{equation}
In the above equation dist$(E,E')$ is the distance between the centroid of the elements $E$ and $E'$, whereas $r_{\textrm{min}}>d_m$ is the size of the filter radius; $d_m$ is the square root of the area of each polygon/element in the discretization $\Tau_h$. Enforcing $r_{\textrm{min}}=1,5d_m$, undesired checkerboard patterns are inhibited; adopting larger values of $r_{\textrm{min}}$ a heuristic  control on the minimum thickness of any member in the design is additionally embedded within the optimization.

\subsection{Compressible elasticity} \label{num:compr}

\begin{figure*}[h]
\centering
\includegraphics[width=7cm]{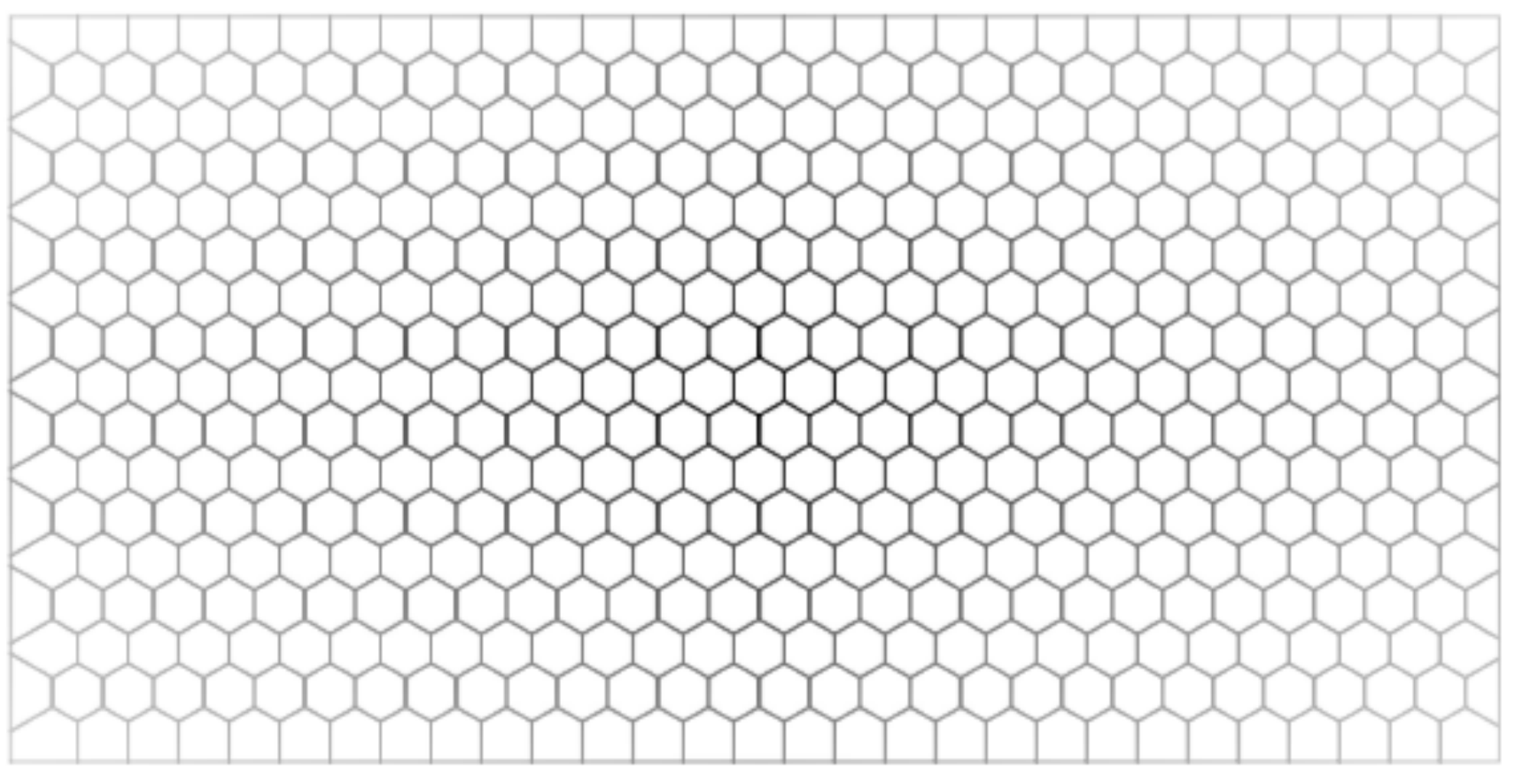}
\includegraphics[width=7cm]{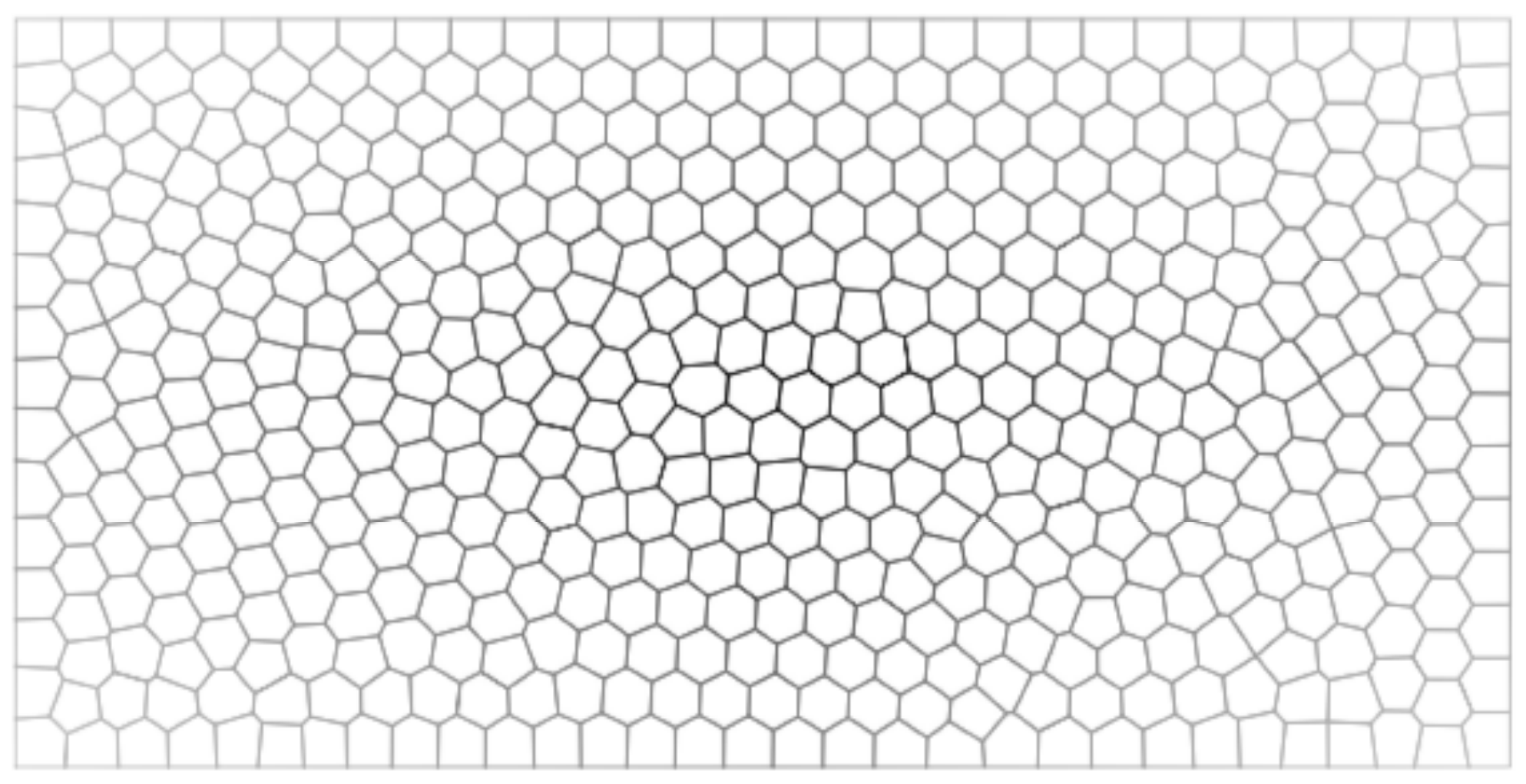}
\caption{Examples of structured (left) and unstructured (right) polygonal grids consisting of $501$ elements.}
\label{fig:grid}
\end{figure*}

In this section  a set of numerical simulations is performed to investigate the features of the proposed VEM--based procedure when addressing topology optimization governed by compressible elasticity. Structured and unstructured polygonal grids have been employed to discretize the design domain. Figure \ref{fig:grid} shows examples of grids including $501$ polygonal elements. The achieved VEM--based results are compared with analogous ones obtained by employing the classical bilinear displacement--based finite elements on Cartesian meshes, see \cite{4}. 
A linear elastic isotropic material is considered in the simulations, assuming Young modulus $E = 1$ and Poisson's ratio $\nu=0.3$. In the whole set of minimizations, the volume fraction of available material is $V_f=0.3$. Different values of filter radius are considered.

\begin{figure*}[h]
\vspace{0.25cm}
\centering\includegraphics[width=12.0cm]{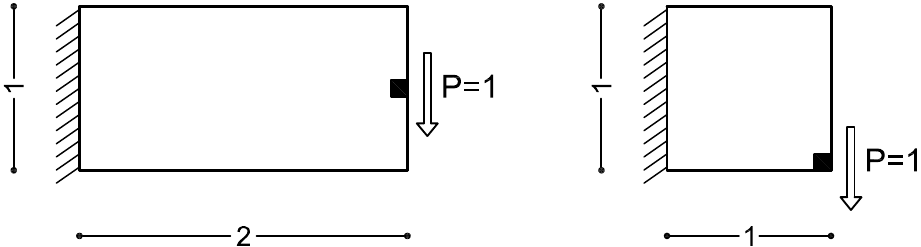}\vspace{0.25cm}
\caption{Compressible elasticity. Geometry and boundary conditions for the numerical simulationsof Example 1 (rectangular cantilever, left) and Example 2 (square cantilever, right). In each case $P=1$ is the intensity of the unitary traction ${\bf f}_t$ applied as a nodal force oriented as indicated by the arrow.}\label{fig:ex}
\end{figure*}

\subsubsection{Example 1: rectangular cantilever}
\begin{figure}[h]
\centering\includegraphics[width=8.0cm]{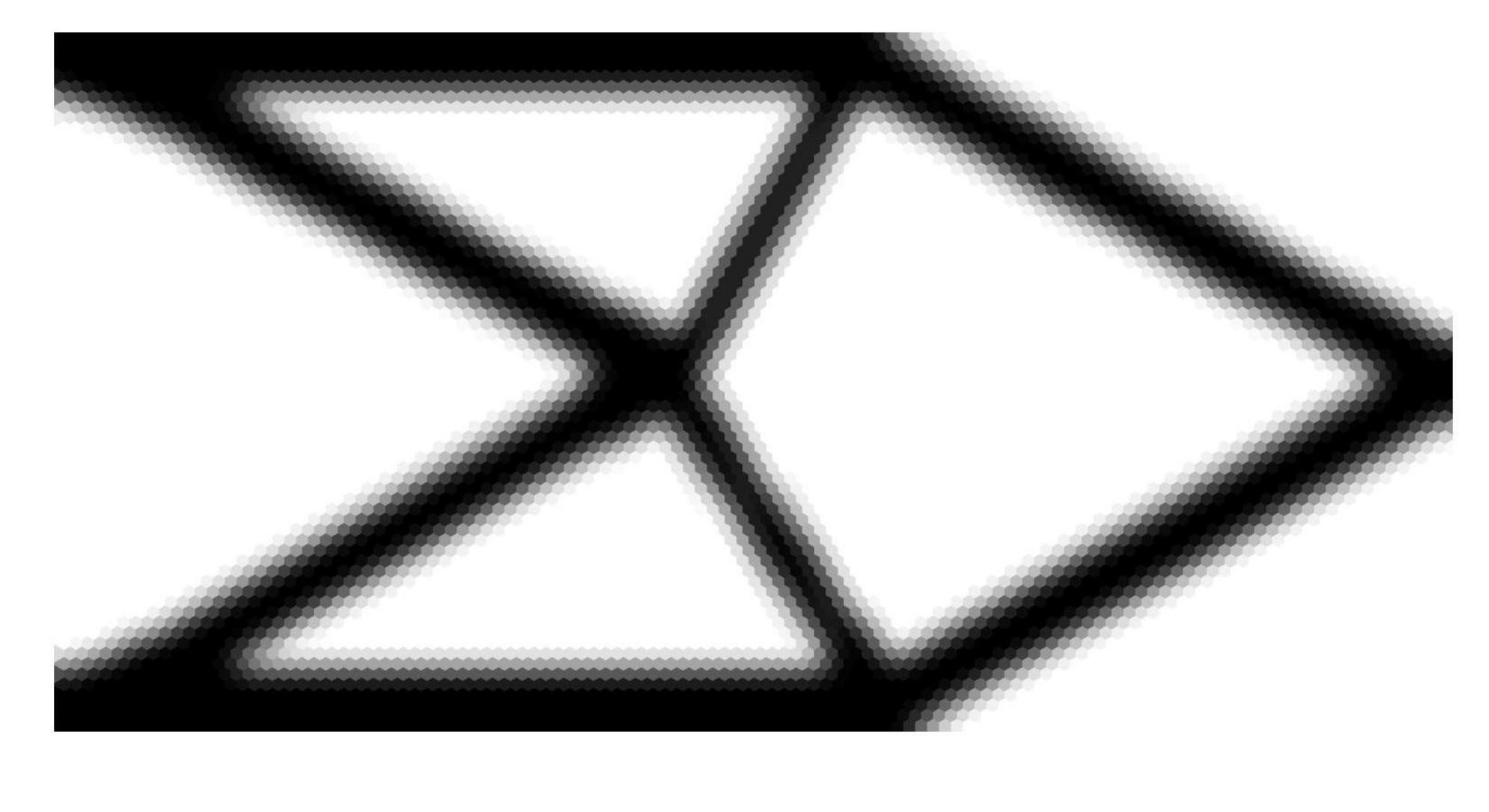}\hspace{-1.25cm}(a)\hspace{0.75cm}\includegraphics[width=8.0cm]{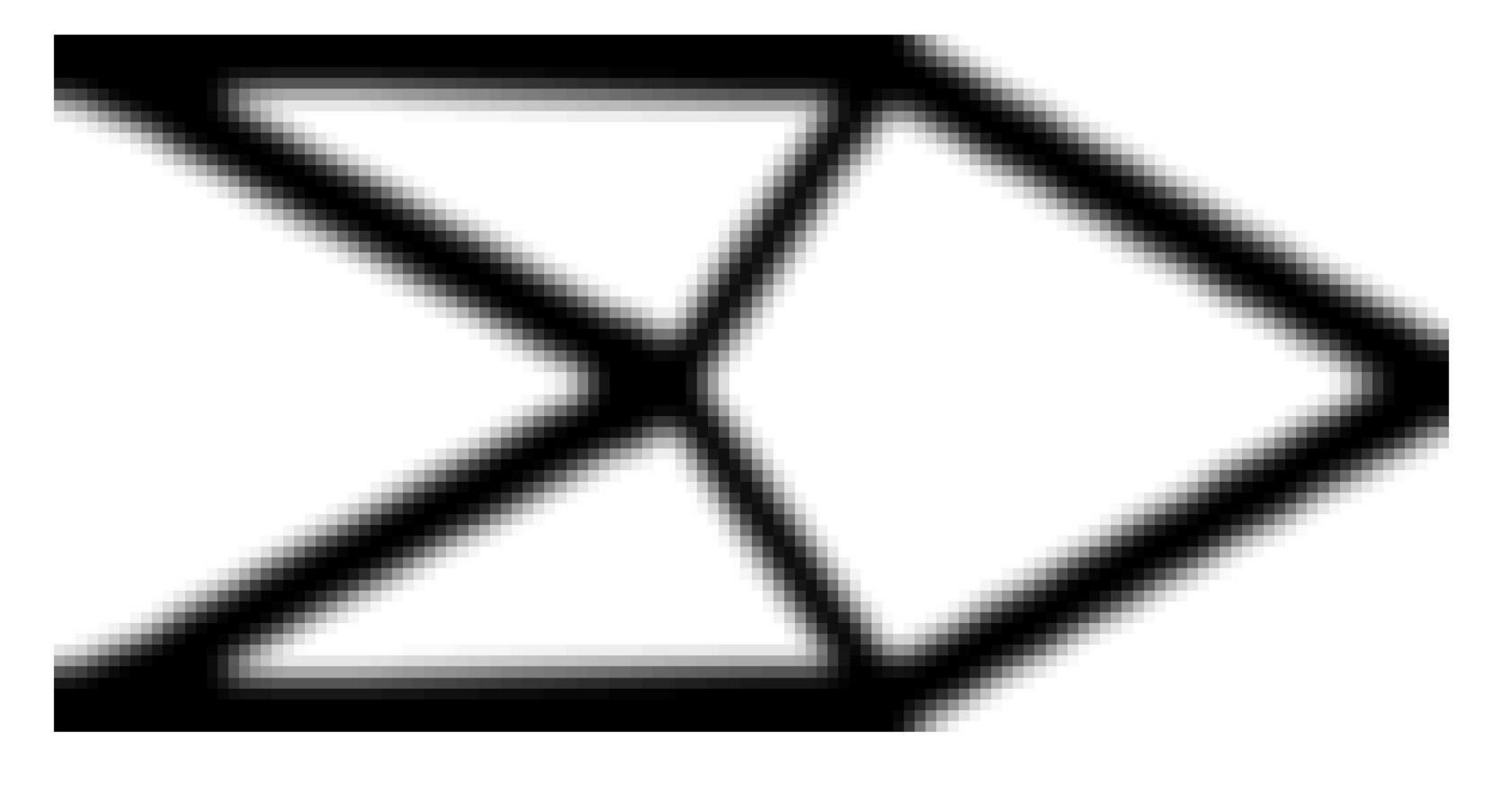}\hspace{-1.25cm}(b)\\
\vspace{0.3cm}
\caption{Compressible elasticity. Example 1: rectangular cantilever. Optimal topologies computed on structured meshes with $2^6$ elements along the thickness of the cantilever: proposed VEM sed formulation (a), bilinear displacement--based formulation (b).} \label{fig:11}
\end{figure}

The first design problem refers to the rectangular cantilever represented in Figure \ref{fig:ex}(left). 
A reference solution is conventionally obtained implementing bilinear displacement--based finite elements on a Cartesian grid consisting of 8192 squares ($2^6$ elements lie along the thickness of the cantilever), along with a filter radius$r_{\textrm{min}}=3.0d_m$. As shown in Figure \ref{fig:11}(b), a truss--like structure arises: inclined members carry shear forces whereas horizontal ones cope with bending actions; both set of members undergo axial stresses. Figure \ref{fig:11}(a) shows the optimal design found through the proposed VEM--based minimization algorithm on a structured mesh consisting of 7990 polygonal elements ($2^6$ elements lie along the thickness of the cantilever). The same filter radius is adopted, e.g. $r_{\textrm{min}}=3.0d_m$. The result obtained with our VEM-based method is substantially equal to the one found by the classical FEM--based approach.\\

\begin{figure}[h]
\centering\includegraphics[width=8.0cm]{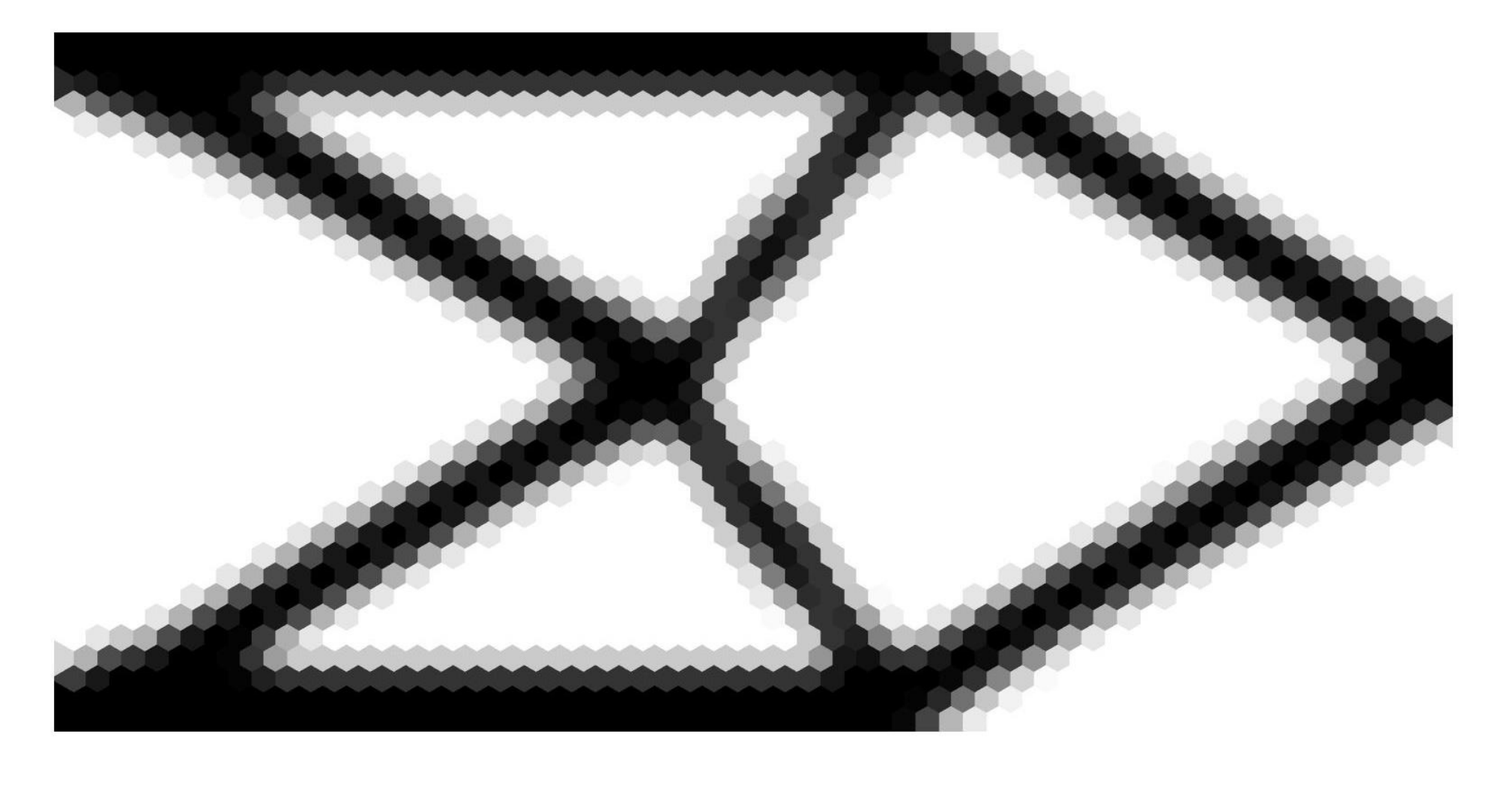}\hspace{-1.25cm}(a)\hspace{0.75cm}\includegraphics[width=8.0cm]{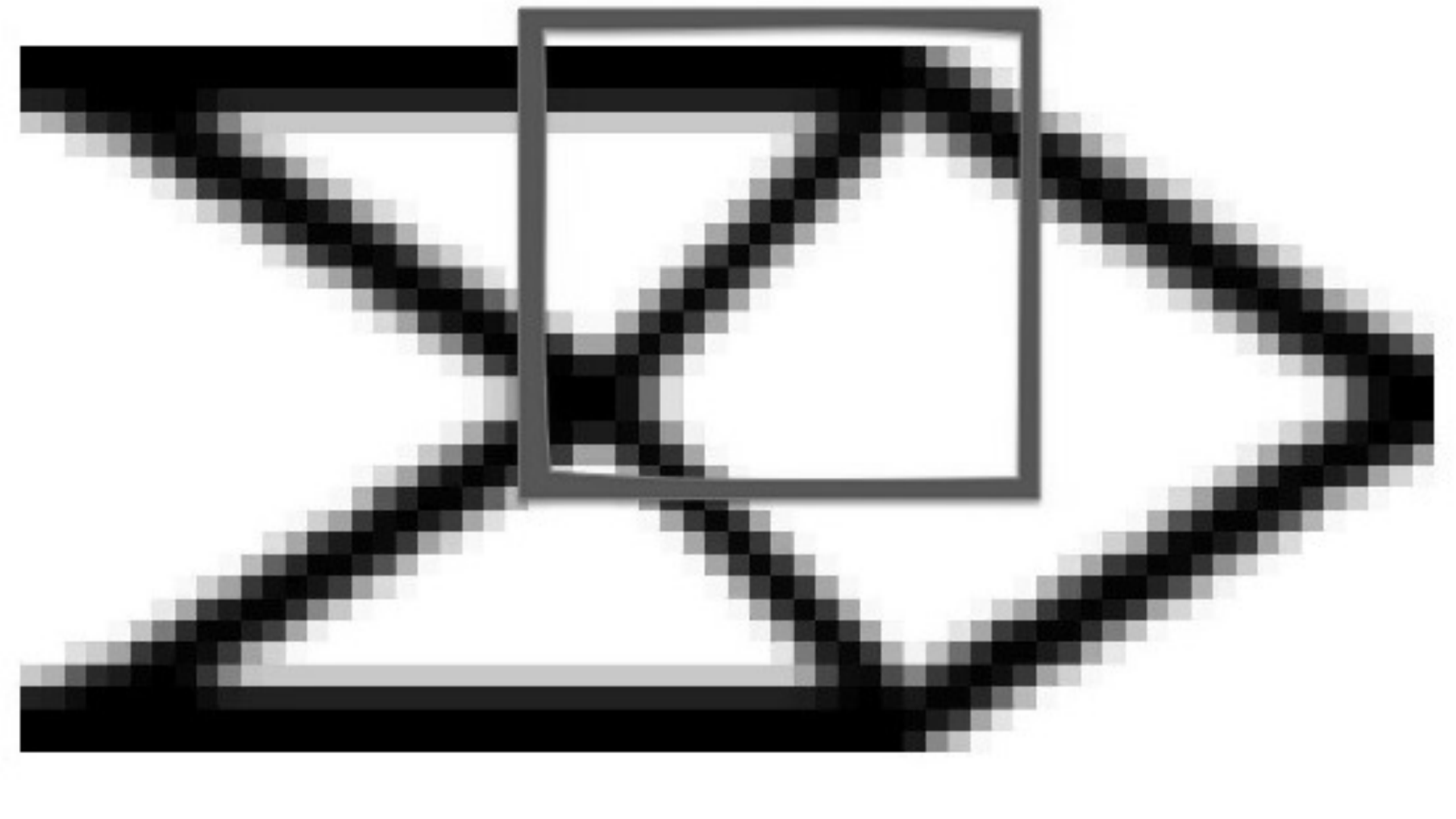}\hspace{-1.25cm}(b)\\
\vspace{0.3cm}
\caption{Compressible elasticity. Example 1: rectangular cantilever. Optimal topologies computed on structured meshes with $2^5$ elements along the thickness of the cantilever: proposed VEM formulation (a), bilinear displacement--based formulation (b).} \label{fig:12}
\end{figure}

\begin{figure}[h]
\centering\includegraphics[width=8.0cm]{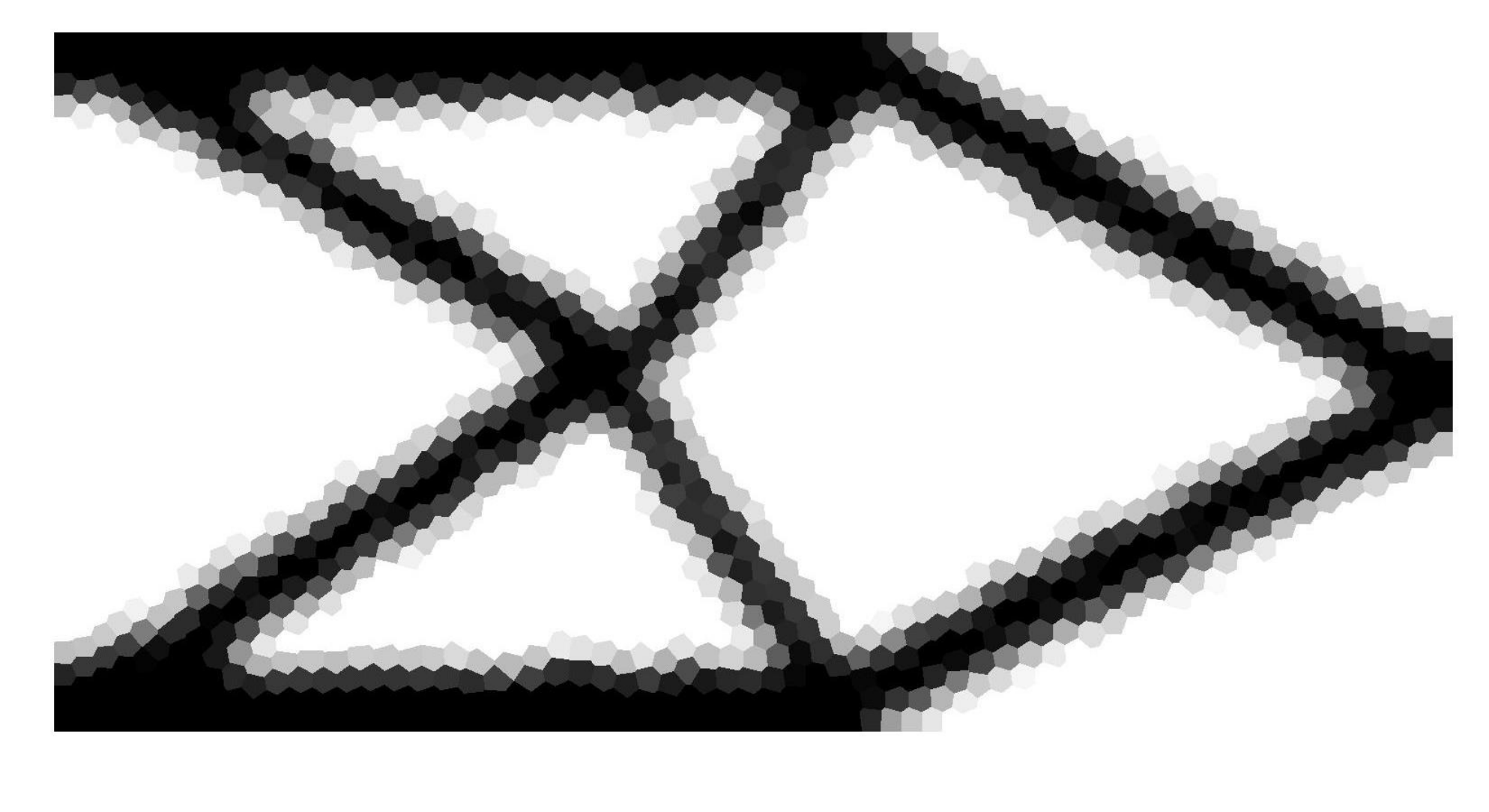}\hspace{-1.25cm}\\
\vspace{0.3cm}
\caption{Compressible elasticity. Example 1: rectangular cantilever. Optimal topology achieved through the proposed VEM formulation for an unstructured mesh with $2^5$ elements along the thickness of the cantilever.} \label{fig:13}
\end{figure}

An additional set of numerical simulations is performed on a coarser discretization adopting $2^5$ elements along the thickness of the cantilever. The same filter radius implemented in the previous investigations is assumed. Figure \ref{fig:12}(a) shows the optimal layout found through the proposed VEM--based approach on a structured mesh consisting of 2006 polygonal elements, whereas Figure \ref{fig:12}(b) refers to the optimal design achieved by the bilinear displacement--based approach on a regular mesh of 2048 square elements. Although the main layout of Figure \ref{fig:11} is recovered in both pictures, the FEM--based design is affected by an unexpected variation in the inclination of the thinner reinforcing braces. Due to the limited amount of available material ($V_f=0.3$) and the rough mesh of square elements, the optimizer gets stuck in a final layout with 45--degree inclination, a mesh--dependent local minimum. The VEM--based result is not affected by such a numerical instability. Indeed, the VEM--based algorithm succeeds in finding the expected layout even in case of coarse unstructured meshes; Figure \ref{fig:13} shows the optimal design obtained on a discretization consisting of 2048 polygonal elements for the same filter radius $r_{\textrm{min}}$ as above.

\begin{figure}[h]
\centering\includegraphics[width=8.0cm]{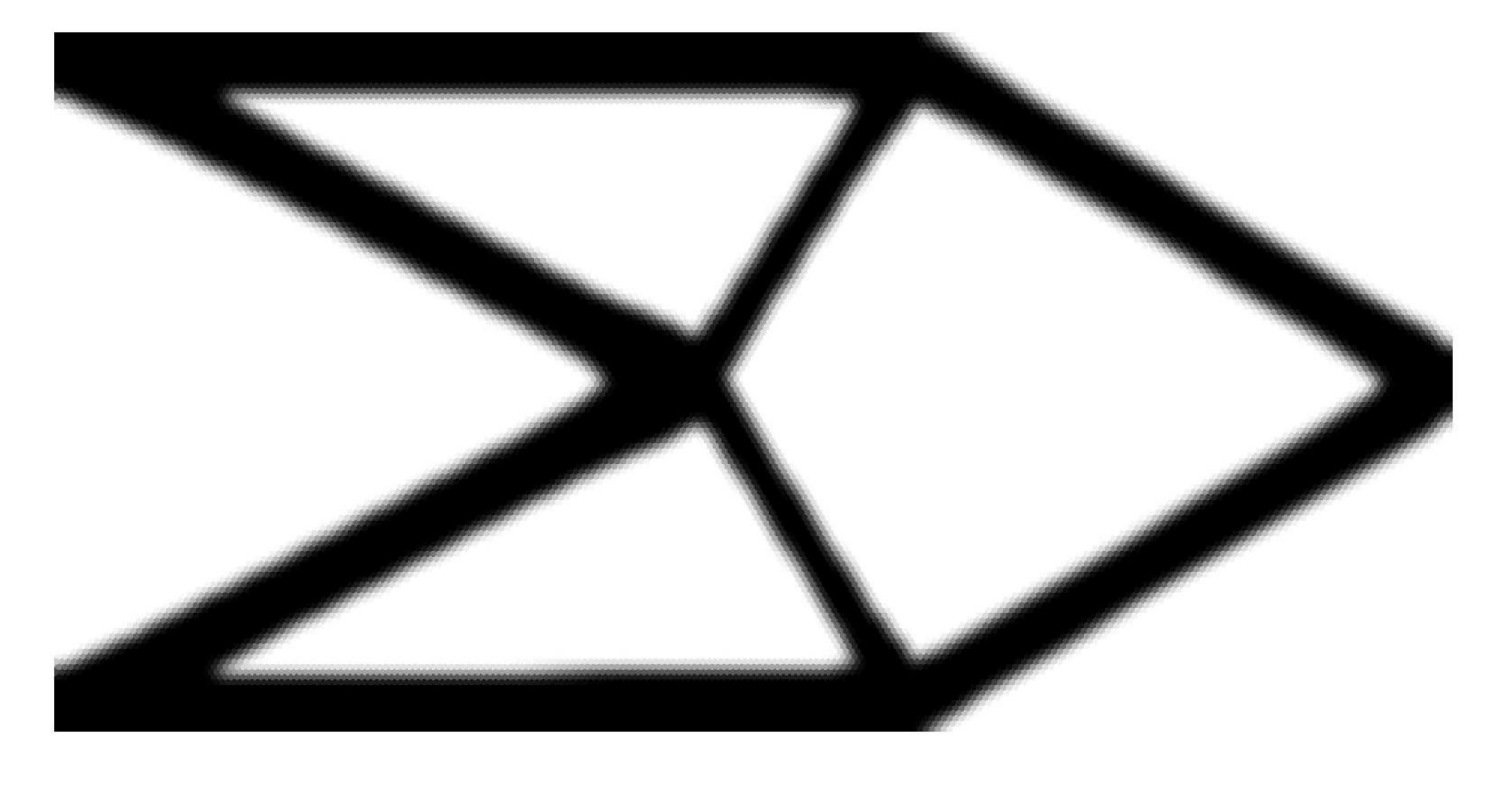}\hspace{-1.25cm}(a)\hspace{0.75cm}\includegraphics[width=8.0cm]{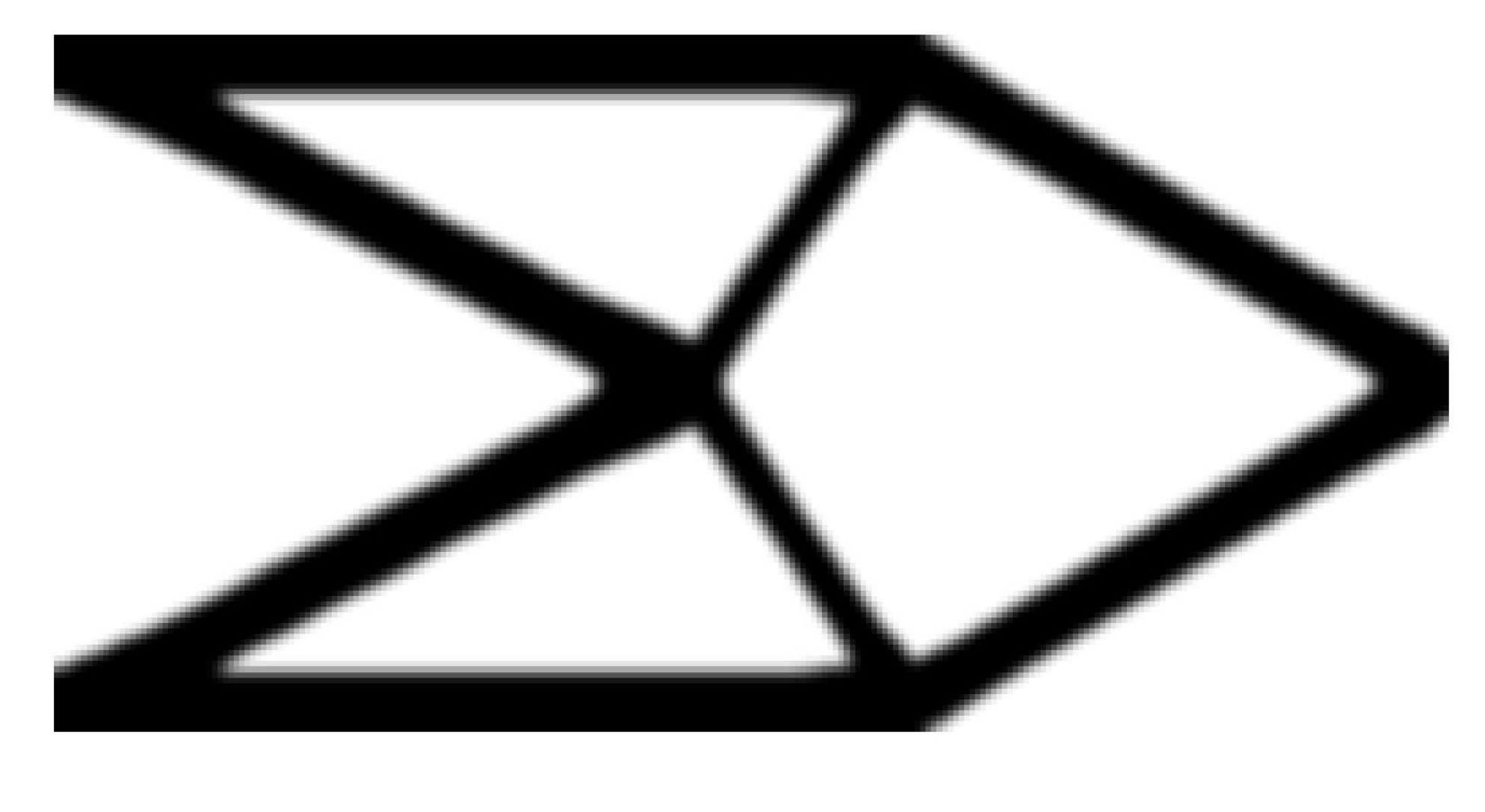}\hspace{-1.25cm}(b)\\
\vspace{0.3cm}
\caption{Compressible elasticity. Example 1: rectangular cantilever. Optimal topologies computed on structured meshes with $2^7$ elements along the thickness of the cantilever: proposed VEM formulation (a), bilinear displacement--based formulation (b).} \label{fig:14}
\end{figure}

Figure~\ref{fig:14} provides a comparison between the VEM--based approach (see Fig.~\ref{fig:14}(a)) and the displacement--based one (see Fig.~\ref{fig:14}(b)) for fine regular meshes of 32028 and 32768 elements, respectively ($2^7$ elements lie along the thickness of the cantilever). Notwithstanding the adopted smaller filter radius $r_{\textrm{min}}=1.5d_m$, the achieved optimal layouts are almost identical.


\subsubsection{Example 2: square cantilever}
The second investigation addresses the square cantilever shown in Figure \ref{fig:ex}(right). An assessment of the VEM--based topology optimization method is provided, adopting unstructured grids of polygonal elements as those shown in Figure \ref{fig:grid}(b).

\begin{figure}[h]
\vspace{0.1cm}
\centering\includegraphics[width=6.0cm]{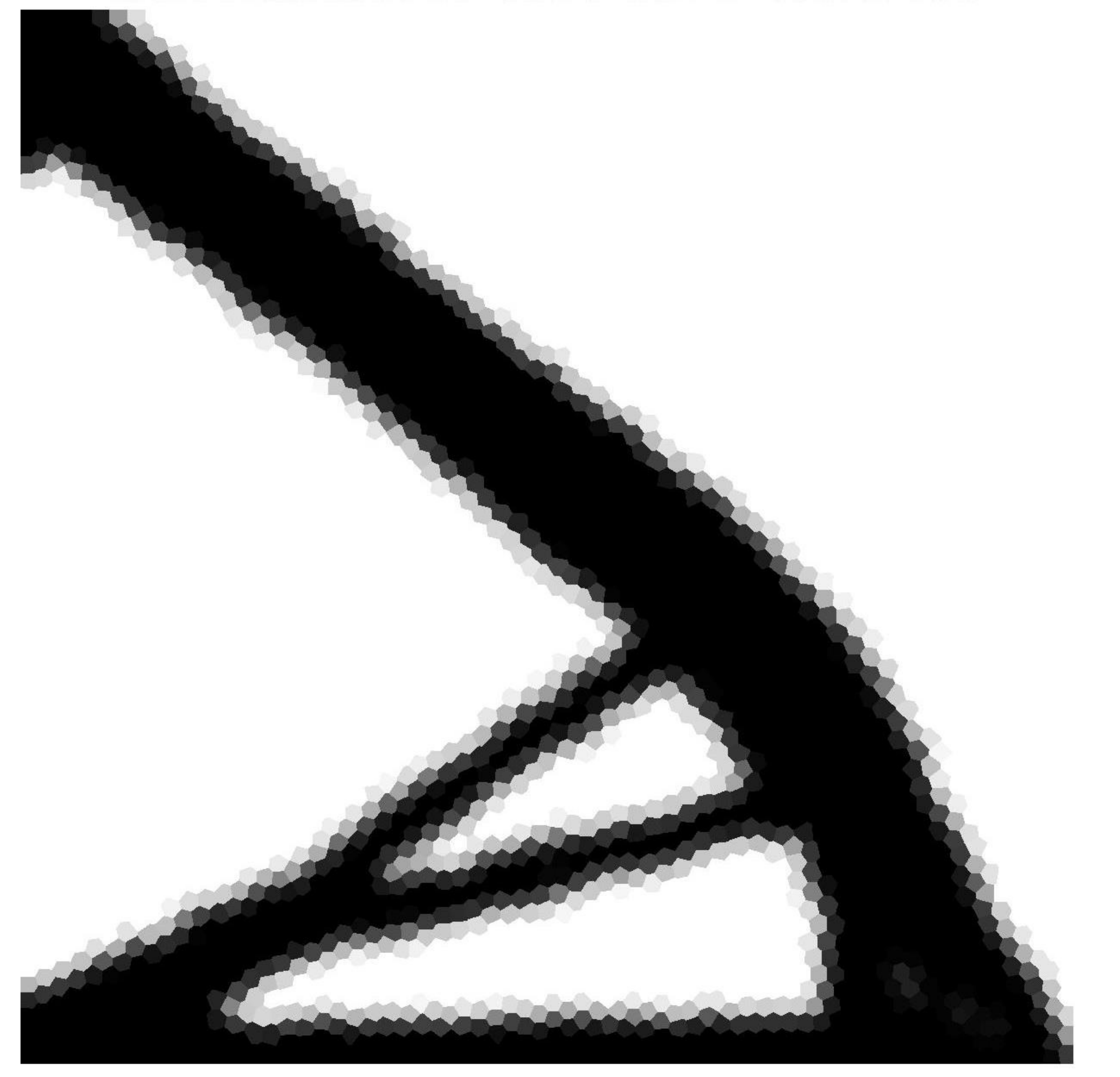}\hspace{-1.25cm}
\\
\vspace{0.1cm}
\caption{Compressible elasticity. Example 2: square cantilever. Optimal topology computed with the proposed VEM formulation on an unstructured mesh with $2^6$ elements along the thickness of the cantilever and filter radius $r_{\textrm{min}}=1.5d_m\approx1.5/2^6$.} \label{fig:21}
\end{figure}

\begin{figure}[h]
\vspace{0.1cm}
\centering\includegraphics[width=6.0cm]{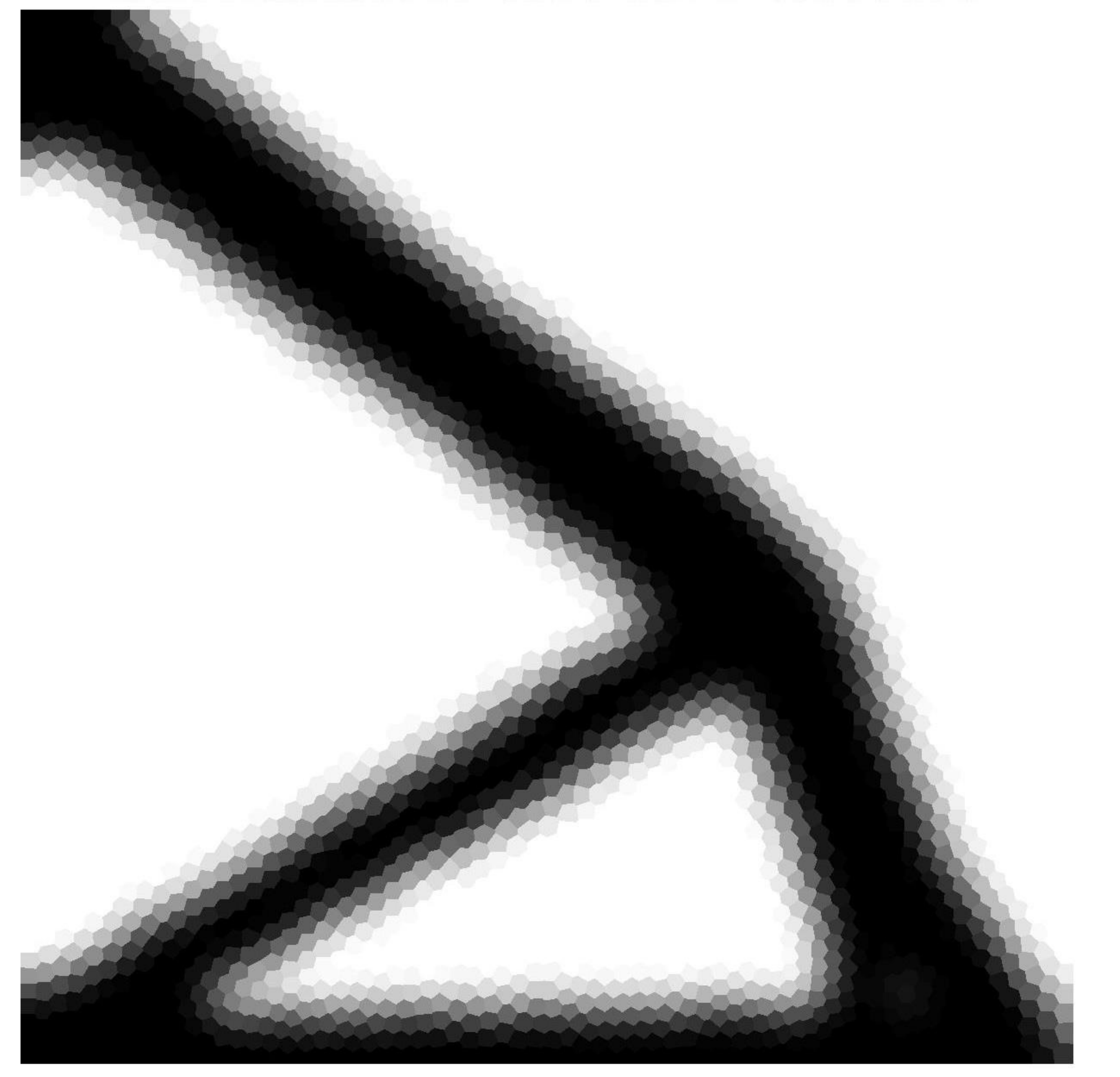}\hspace{-1.25cm}\\
\vspace{0.1cm}
\caption{Compressible elasticity. Example 2: square cantilever. Optimal topology computed with the proposed VEM formulation on an unstructured mesh with $2^6$ elements along the thickness of the cantilever and filter radius $r_{\textrm{min}}=3.0d_m\approx1.5/2^5$.} \label{fig:22}
\end{figure}

First, an unstructured discretization accounting for 4096 elements ($2^6$ along the thickness of the cantilever) is used. The adopted filter radius reads $r_{\textrm{min}}=1.5d_m$, being $d_m$, as before, the square root of the average area of the polygonal elements in the unstructured grid. Figure \ref{fig:21} shows the achieved optimal layout, a truss--like structure whose central node receives two thick tensile--stressed trusses along with two thin compressive--stressed bars.
Figure \ref{fig:22} shows the optimal solutions achieved for an increased value of the enforced filter radius, i.e. $r_{\textrm{min}}=3.0d_m$. A simpler design arises that is made of two ties and one big strut, in full agreement with the well--known solution of this benchmark problem, see e.g. \cite{4}.\\

\begin{figure}[h]
\vspace{0.1cm}
\centering\includegraphics[width=6.0cm]{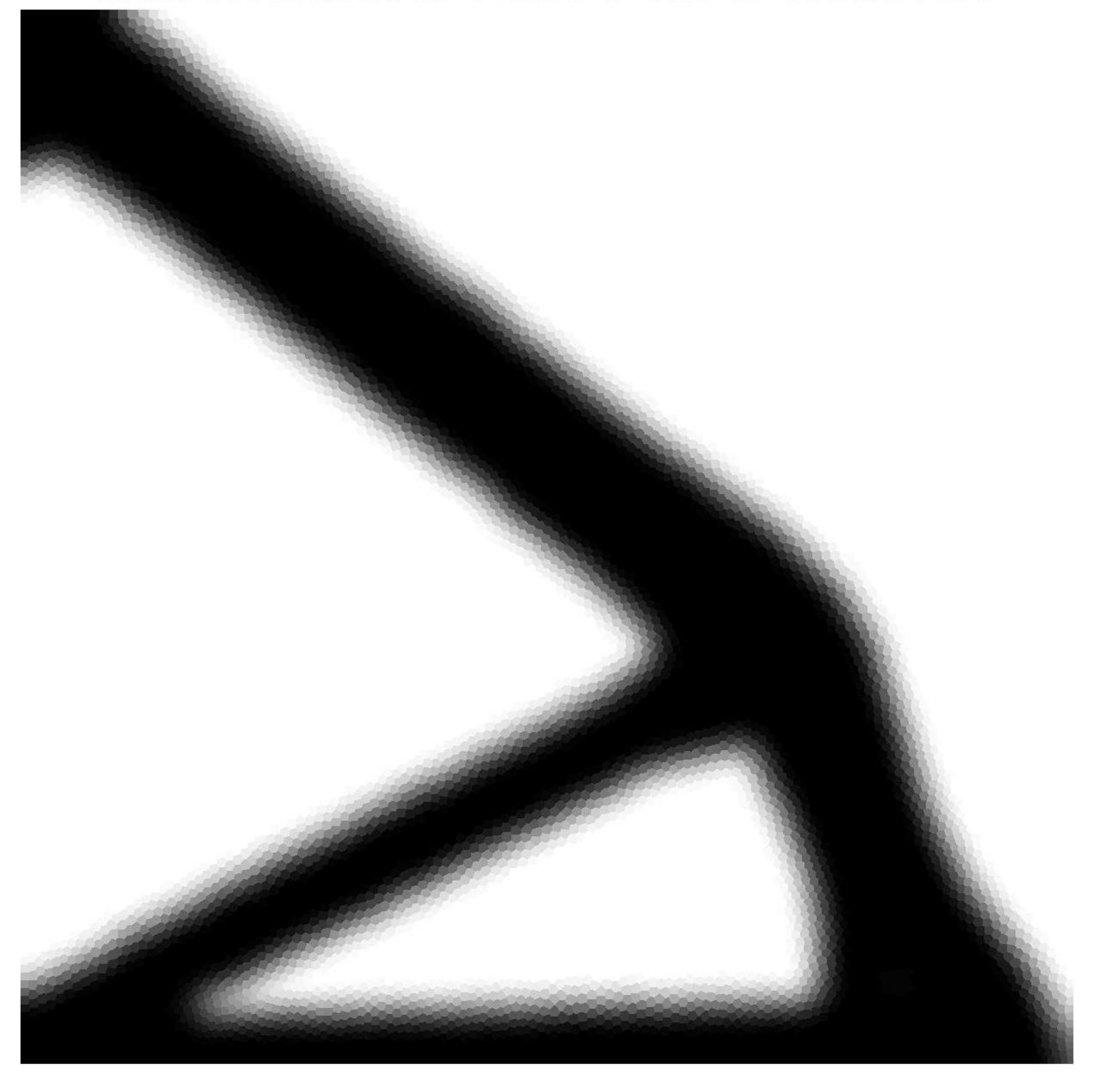}\hspace{-1.25cm}\\
\vspace{0.1cm}
\caption{Compressible elasticity. Example 2: square cantilever. Optimal topology computed with the proposed VEM formulation  on an unstructured mesh with $2^7$ elements along the thickness of the cantilever and filter radius $r_{\textrm{min}}=6.0d_m\approx1.5/2^5$.} \label{fig:23}
\end{figure}

Finally, Figure \ref{fig:23} provides the optimal layout computed when a finer discretization with $2^7$ elements along the thickness of the cantilever is implemented. The overall number of polygonal elements is 16384. The filter radius is $r_{\textrm{min}}=6.0d_m$, which is nearly the value used for the result shown in Figure \ref{fig:22}. No mesh dependence affects the proposed VEM--based formulation, since the same optimal layout is found in both figures.


\subsubsection{Example 3: circle loaded with four point load} 
\begin{figure}
\vspace{0.25cm}
\centering\includegraphics[width=0.47\textwidth]{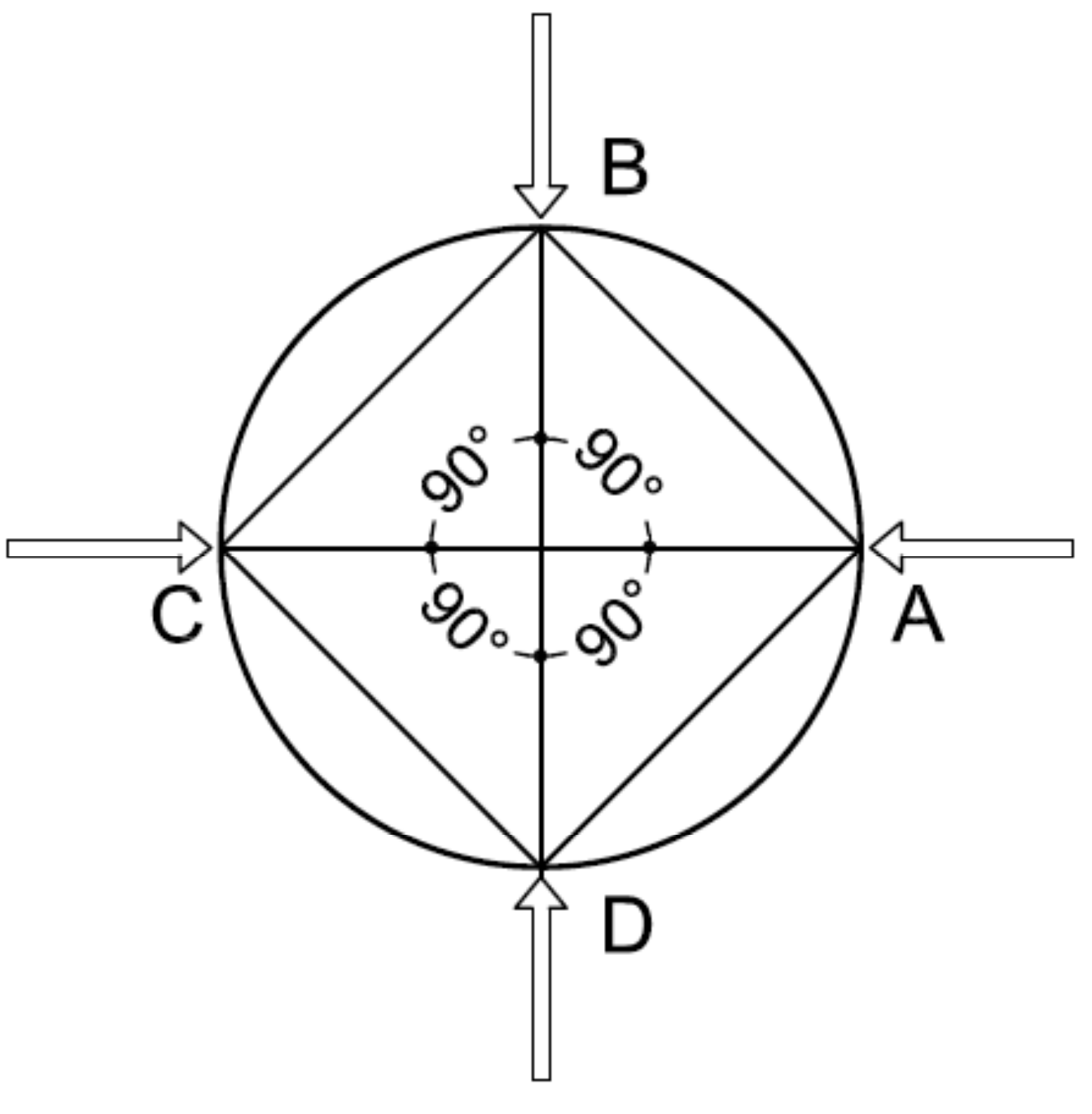}\vspace{0.25cm}

\caption{Compressible elasticity. Geometry and boundary conditions for the numerical simulations of 
Example 3 (circle loaded with four point load).}\label{fig:ex:2}
\end{figure}

Let us consider a circular lamina that is loaded by a set of self--balanced forces applied at points A,B,C,D as shown in Figure~ \ref{fig:ex:2}.\\

\begin{figure}
\centering
\includegraphics[width=0.47\textwidth]{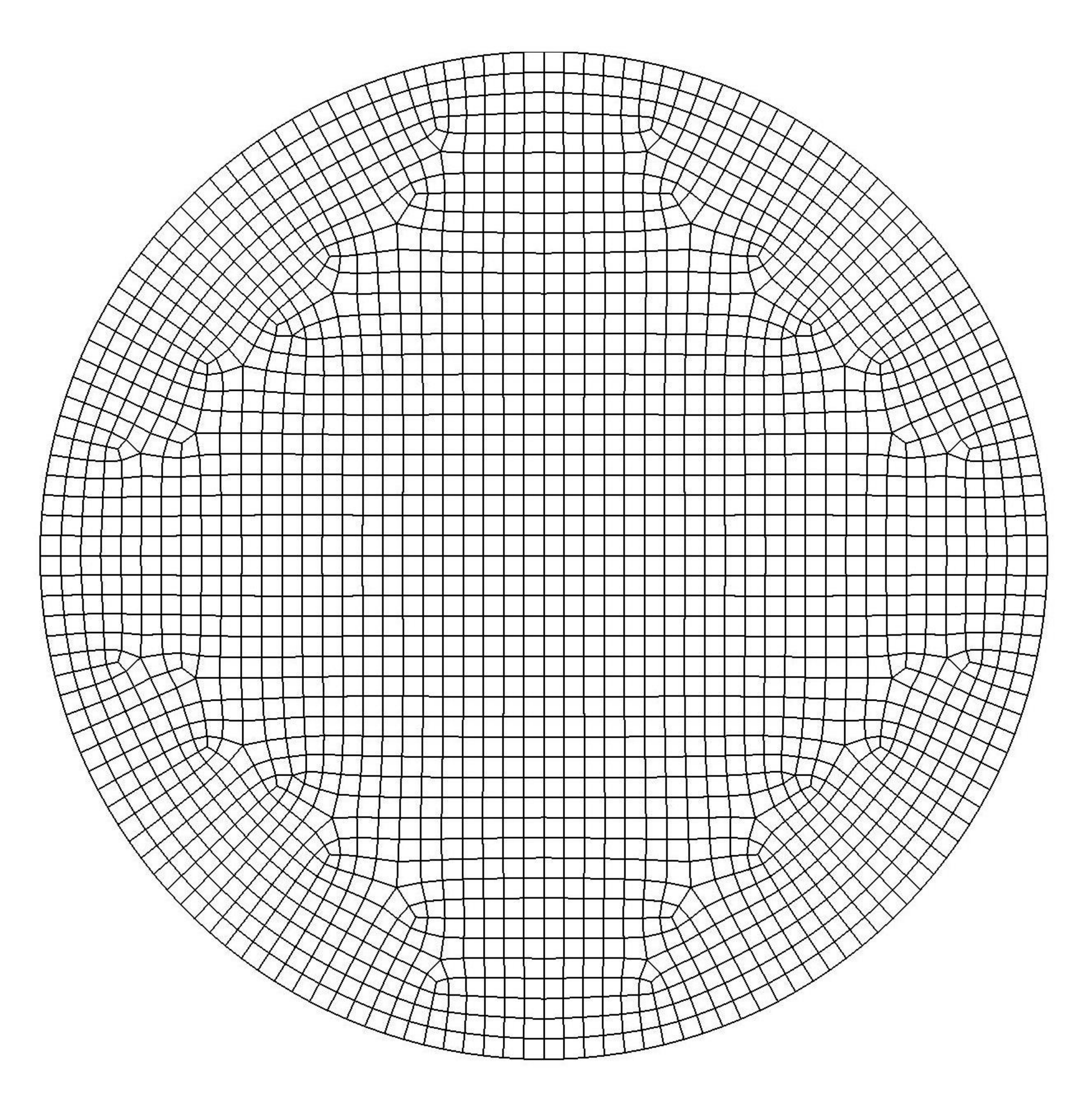}
\caption{Compressible elasticity. Example 3: circle loaded with four point load. Mesh of $2168$ quadrilateral elements built with commercial code {\tt Strand7} }\label{Fig:circ:1}
\end{figure}

\begin{figure}
\centering\includegraphics[width=0.47\textwidth]{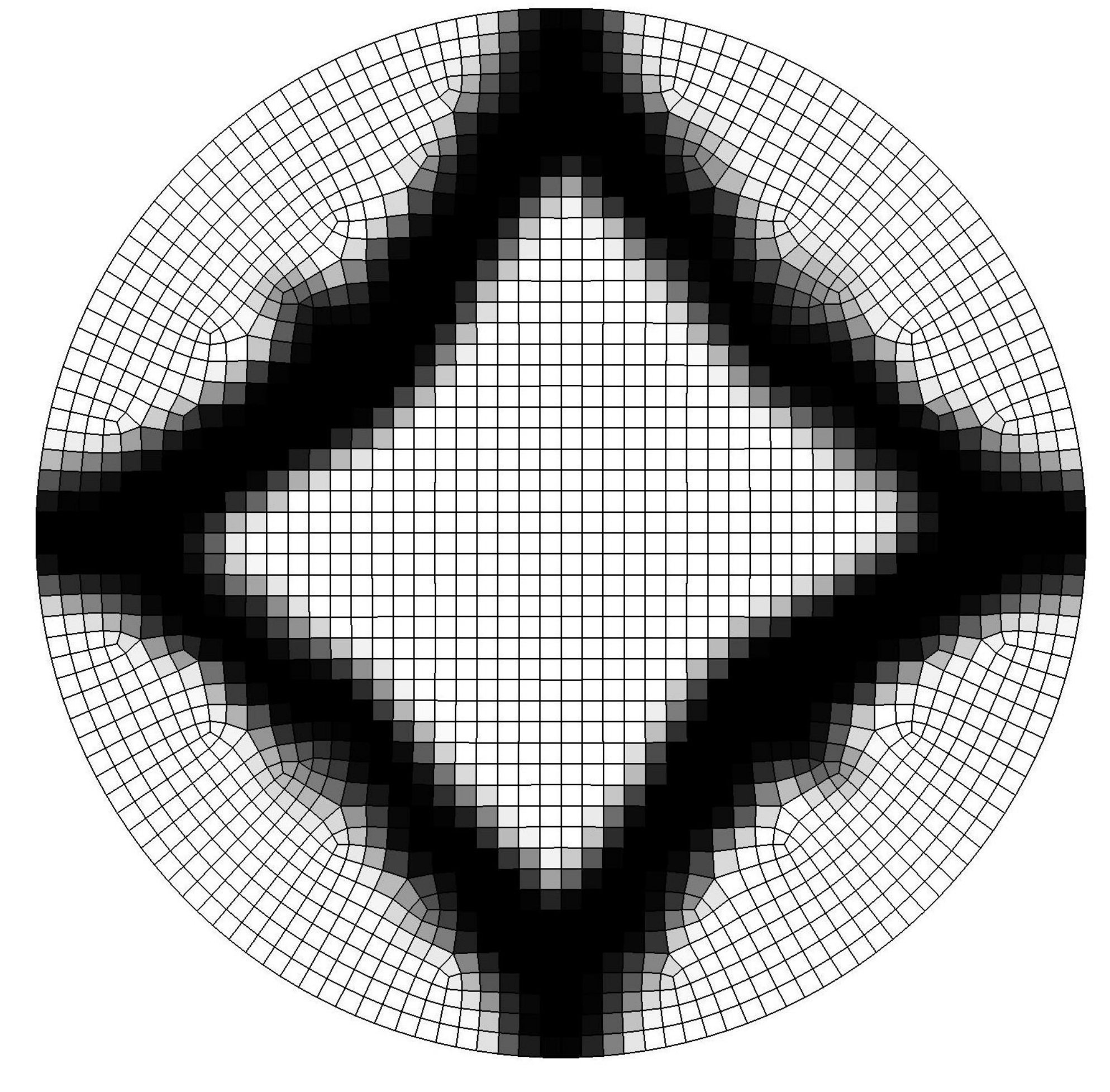}\hspace{-1.25cm}(a)\hspace{0.75cm}%
{\centering\includegraphics[width=0.47\textwidth]{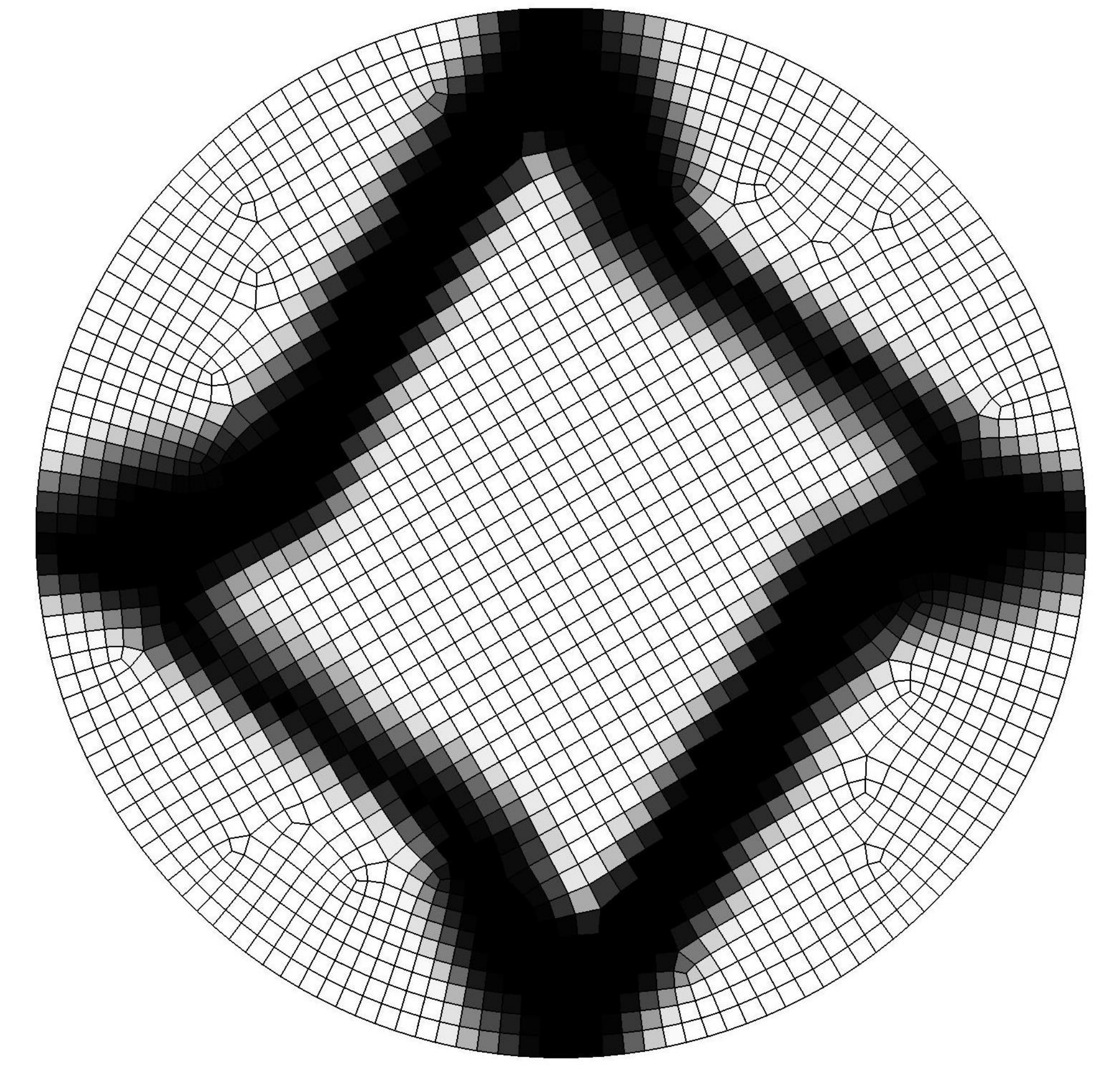}(b)}
\vspace{0.3cm}
\caption{Compressible elasticity. Example 3: circle loaded with four point load. Final configuration with bilinear displacement--based formulation: (a) no mesh rotation;  {{(b) mesh rotation of $30$ degrees.}}}\label{Fig:circ:2}

\end{figure}

\begin{figure}
\centering
\includegraphics[width=0.66\textwidth]{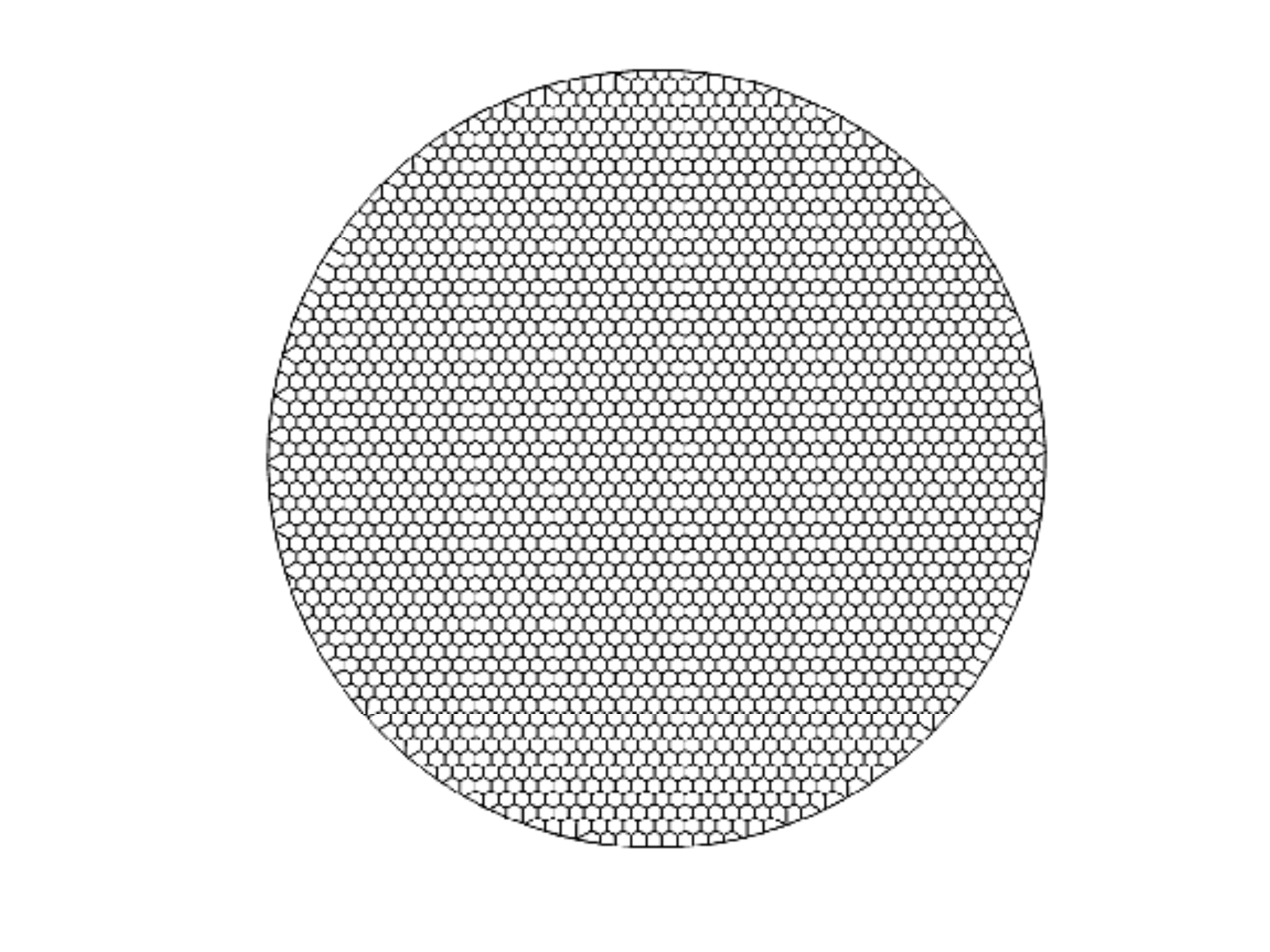}
\caption{Compressible elasticity. Example 3: circle loaded with four point load. Polygonal grid consisting of $2268$ elements built with {\tt Polymesher} \cite{PolyMesher}.}\label{Fig:circ:3}
\end{figure}

The geometry is discretized using the commercial code {\tt Strand7} \cite{strand7} to achieve a mesh of $2168$ quadrilateral elements, see Figure~ \ref{Fig:circ:1}. The conventional FEM--based formulation is adopted to find the volume--constrained minimum compliance solution. First, the load case presented in Figure~ \ref{fig:ex:2} is applied to the mesh in Figure~ \ref{Fig:circ:1}. Then, the same load case is applied to the discretization achieved after a $30$--degree anticlockwise rotation of the original mesh around its centroid. Figure~\ref{Fig:circ:2} shows the result of the topology optimization procedure for the unrotated mesh of standard bilinear finite elements (see Fig.~\ref{Fig:circ:2}(a)) and the rotated one (see Fig.~\ref{Fig:circ:2}(b)). The achieved layouts are remarkably different and point out a lack of robustness of the solution with respect to the considered rotation of the mesh. In both cases a mesh--dependent solution is found that is not a truss--like structure. Due to the geometrical features of the adopted discretizations, non--straight members arise, curved beams in Figure~ \ref{Fig:circ:3}(a) or piecewise linear elements in Figure~ \ref{Fig:circ:3}(b). Both kinds of elements have to cope with bending stresses, meaning that a sub--optimal performance is achieved with respect to any stiff truss--like structure.\\   

\begin{figure}
\centering\includegraphics[width=0.50\textwidth]{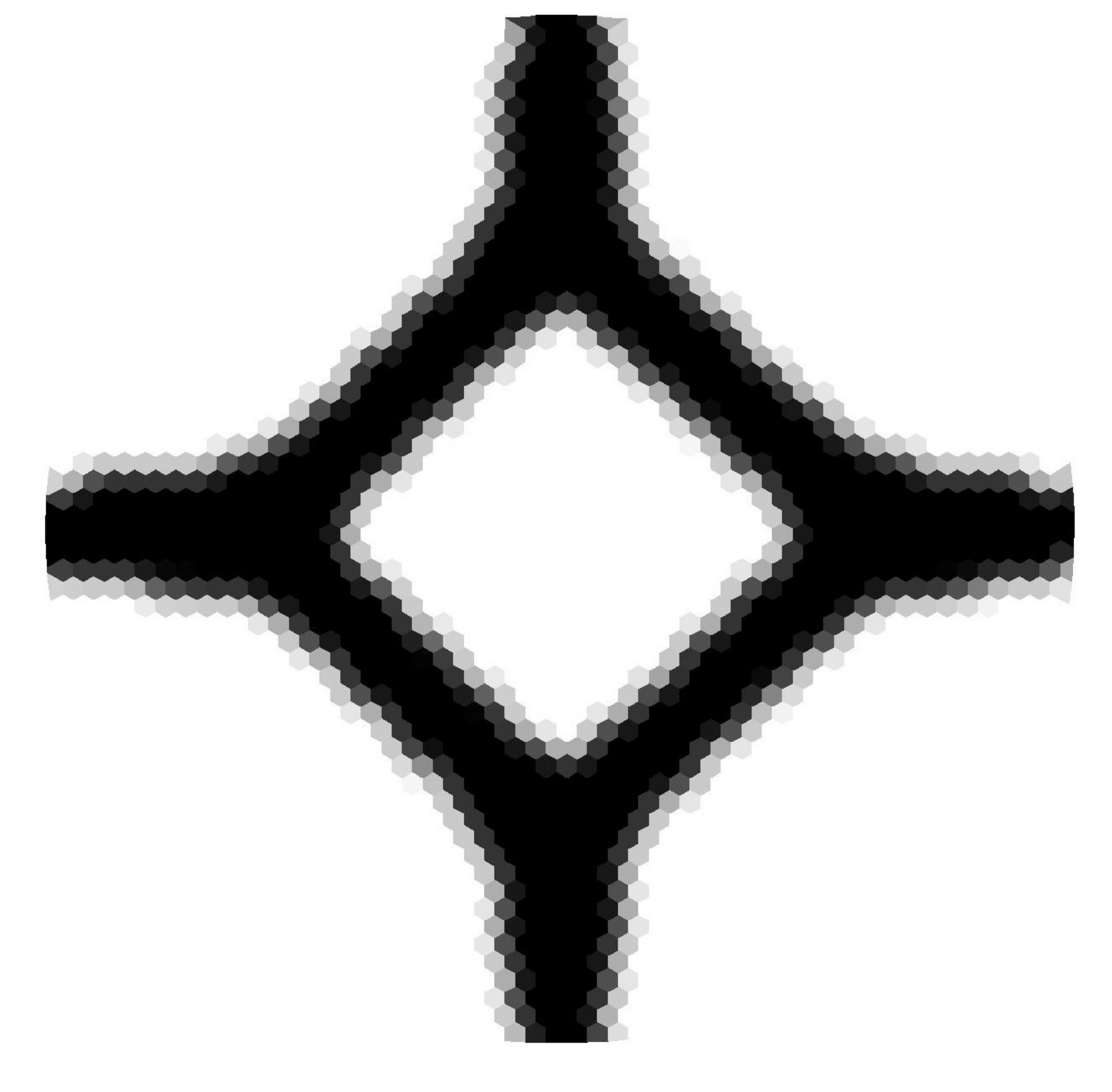}\hspace{-1.25cm}(a){\hspace{0.75cm}\includegraphics[width=0.50\textwidth]{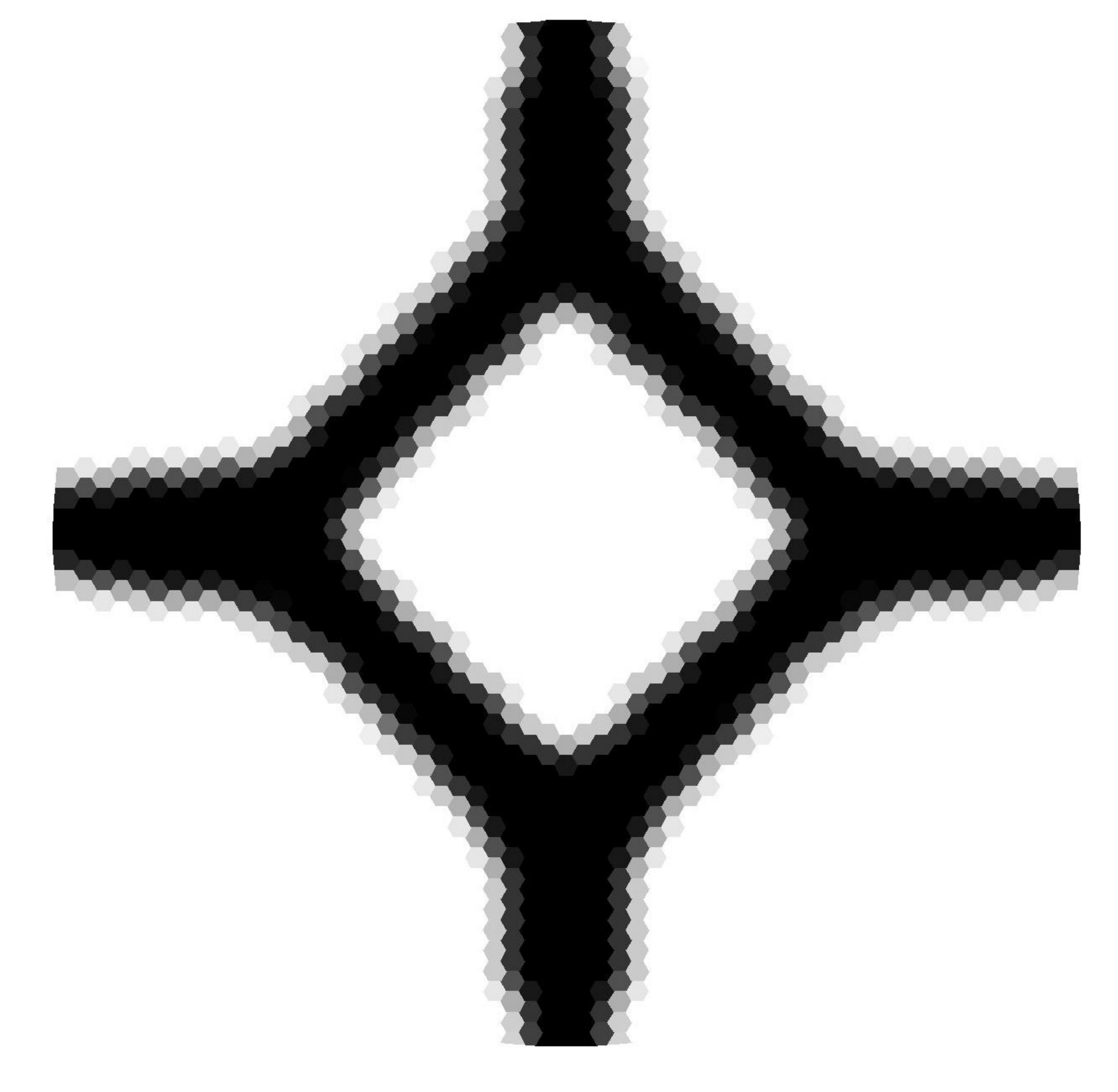}\hspace{-1.25cm}(b)}\\
\vspace{0.3cm}
\caption{Compressible elasticity. Example 3: circle loaded with four point load. Final configuration with VEM on structured meshes with $2268$ elements: (a) no mesh rotation; 
{(b) mesh rotation of $30$ degrees.}} \label{Fig:circ:4}
\end{figure}

A similar investigation is performed adopting the proposed VEM--based approach to solve the considered problem of optimal design. The geometry is discretized using the academic code {\tt Polymesher} \cite{PolyMesher} to achieve a mesh of $2268$ polygonal elements, see Figure~ \ref{Fig:circ:3}. Figure~ \ref{Fig:circ:4} shows the results of the VEM--based topology optimization for the unrotated mesh of Figure~ \ref{Fig:circ:3} (a) and for a discretization achieved after a $30$--degree anticlockwise rotation of the original mesh around its centroid (b). The achieved layouts are very similar to each other: the considered rotation of the mesh induces only minor effects on the solution. In both cases a stiff truss--like structure arises, thus assessing the robustness of the proposed algorithm with respect to mesh rotations.



\subsection{Nearly-incompressible elasticity} \label{num:incompr}
\begin{figure}
\centering
\includegraphics[width=0.45\textwidth]{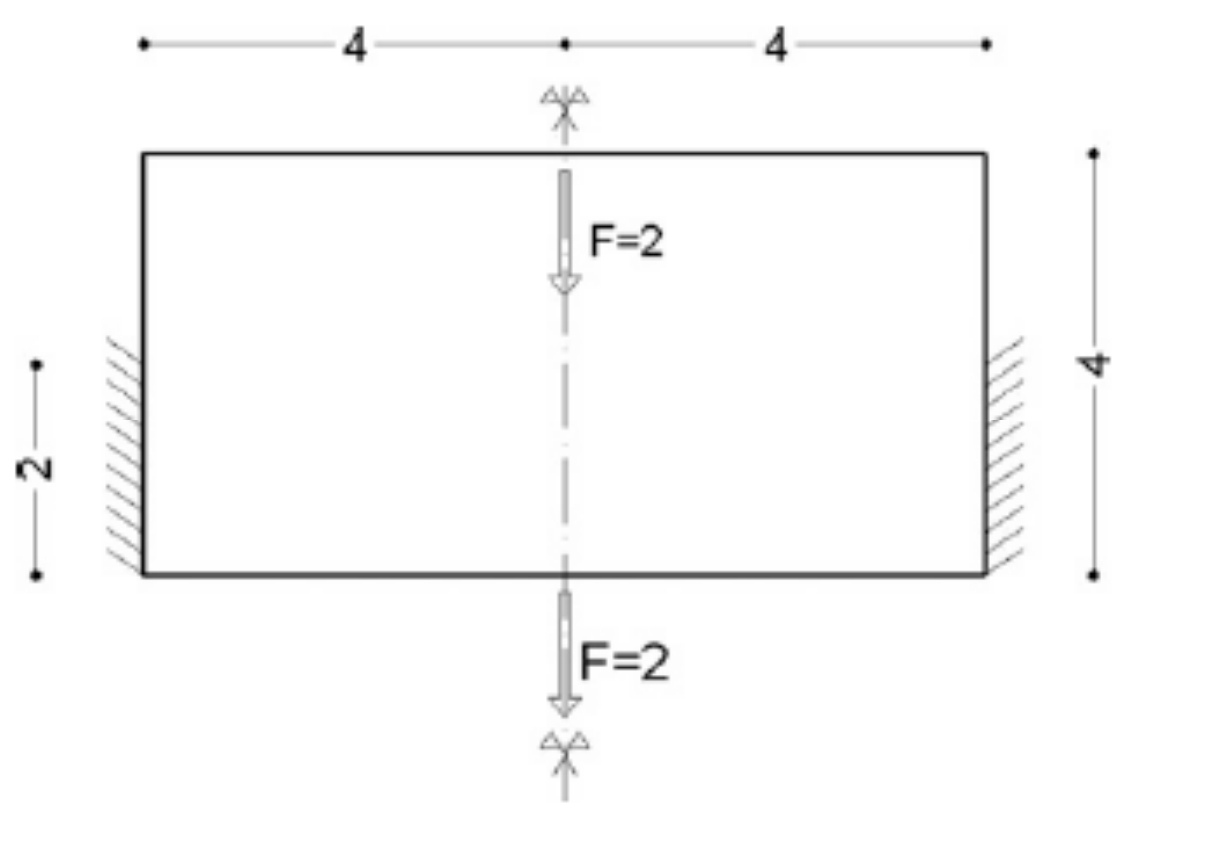}
\includegraphics[width=0.45\textwidth]{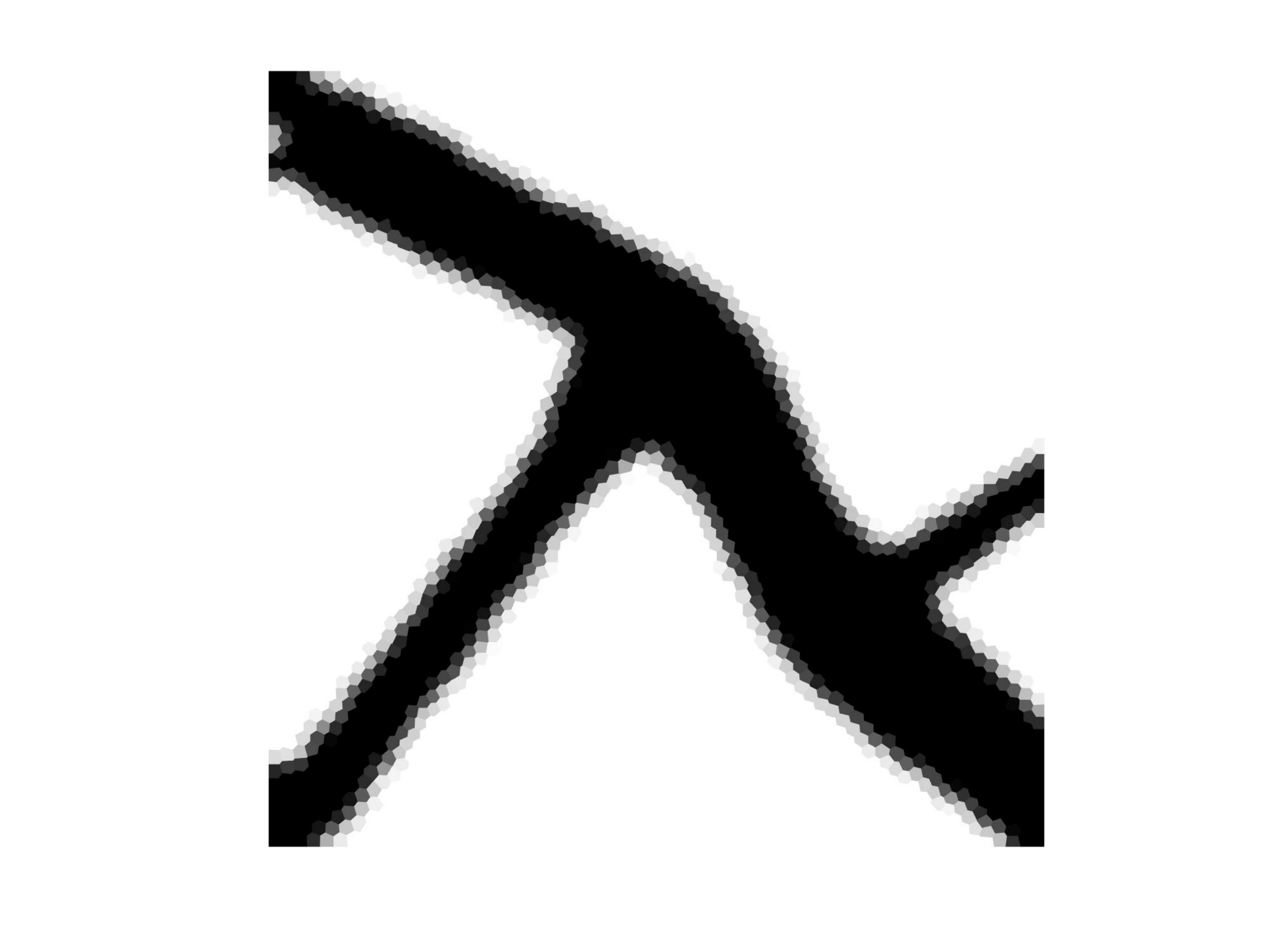}
\caption{Nearly-incompressible elasticity. Left: design domain. Right: optimal design (half by symmetry) on polygonal grid consisting of $4096$ hexagons ($\nu=0.4999999$). }\label{Fig:inc:2}
\end{figure}

Differently from displacement--based finite elements, the adopted VEM approximation is well--suited to cope with the analysis of quasi--incompressible media, since, at least numerically (see Remark \ref{rem:inf-sup}), it satisfies the classical inf-sup stability condition and no locking is expected when dealing with problems assuming plane strain. As originally investigated in \cite{7}, a conventional SIMP--law that uses the same exponent $p$ to approximate the dependence of the modulus $\lambda$ and $\mu$ with respect to the density $\rho$ can fail when addressing the optimal design of (nearly-)incompressible media: region with low density but high stiffness may arise in the solution, thus achieving optimal layouts that are unfeasible from a physical point of view. This problem can be simply overcome adopting a larger penalization on $\lambda$ than $\mu$, e.g. enforcing  $p_{\lambda}=6$ along with $p_{\mu}=3$ instead of $p_{\lambda}=p_{\mu}=3$, see \cite{Bruggi-Venini:07}.\\

Figure~\ref{Fig:inc:2}(left) shows a benchmark problem for nearly--incompressible two--dimensional bodies undergoing plane strain conditions. Figure~\ref{Fig:inc:2}(right) shows the optimal design obtained by our VEM--based procedure (for symmetry reasons only half of the domain has been tackled in the optimization). The achieved layout is in full agreement with those found in the literature, see in particular the results obtained by adopting robust truly--mixed discretizations based on triangles \cite{Bruggi-Venini:07} or square elements \cite{Bruggi:16}.



\subsection{Optimal Stokes flow} \label{num:stokes}
In this section we consider the numerical solution of the discrete problem \eqref{eq:opti3_h} related to the  optimization of Stokes flows.  In particular, we will deal with some classical benchmark examples first introduced in \cite{Borrvall-Petersson:03} together with some suitable variants aiming at highlighting the virtues of the VEM-based topology optimization on polygonal meshes. 

In the sequel, following  \cite{Borrvall-Petersson:03} we will work under the following choice of the penalty function $\alpha(\dens)$ (cf. \eqref{stokes:form}), namely 
\[
\alpha(\dens)=\bar{\alpha} + (\underline{\alpha} - \bar{\alpha} )\dens \frac{1+q}{\dens+q},
\] 
where $q=0.1$, $\underline{\alpha}= 2.5 \mu_0/100^2$ and 
$\bar{\alpha} =  2.5 \mu_0/0.01^2$.

The profile of the prescribed non-zero velocity at the boundary is parabolic, i.e. the magnitude of the velocity can be written as $g^*(1-(2s/l)^2)$, with $s\in [-l/2,l/2]$ where $g^*$ is the maximum value, whereas $l$ is the length of the boundary part where the parabolic profile is prescribed. 

Throughout all the following numerical experiments we choose $\mu_0=1$ and $\lambda_0=10^3$.

\subsubsection{Optimal pipe} 
\begin{figure}
\centering
\includegraphics[width=0.45\textwidth]{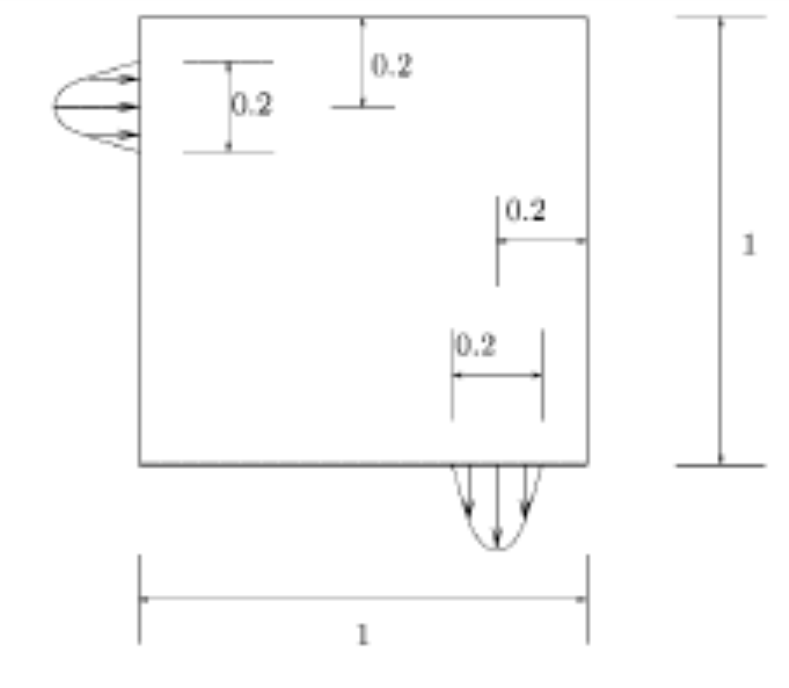}\\
\includegraphics[width=0.45\textwidth]{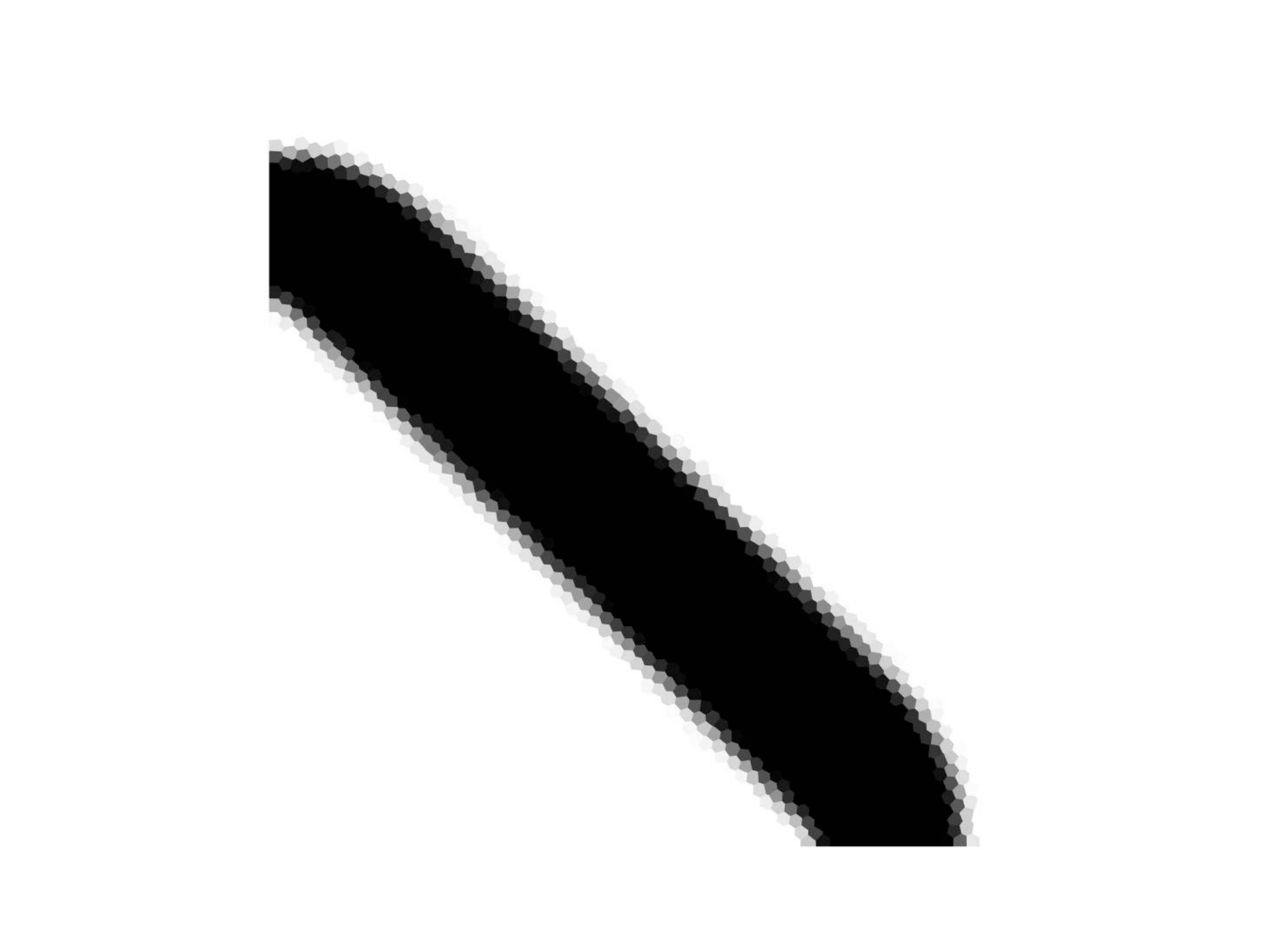}
\includegraphics[width=0.45\textwidth]{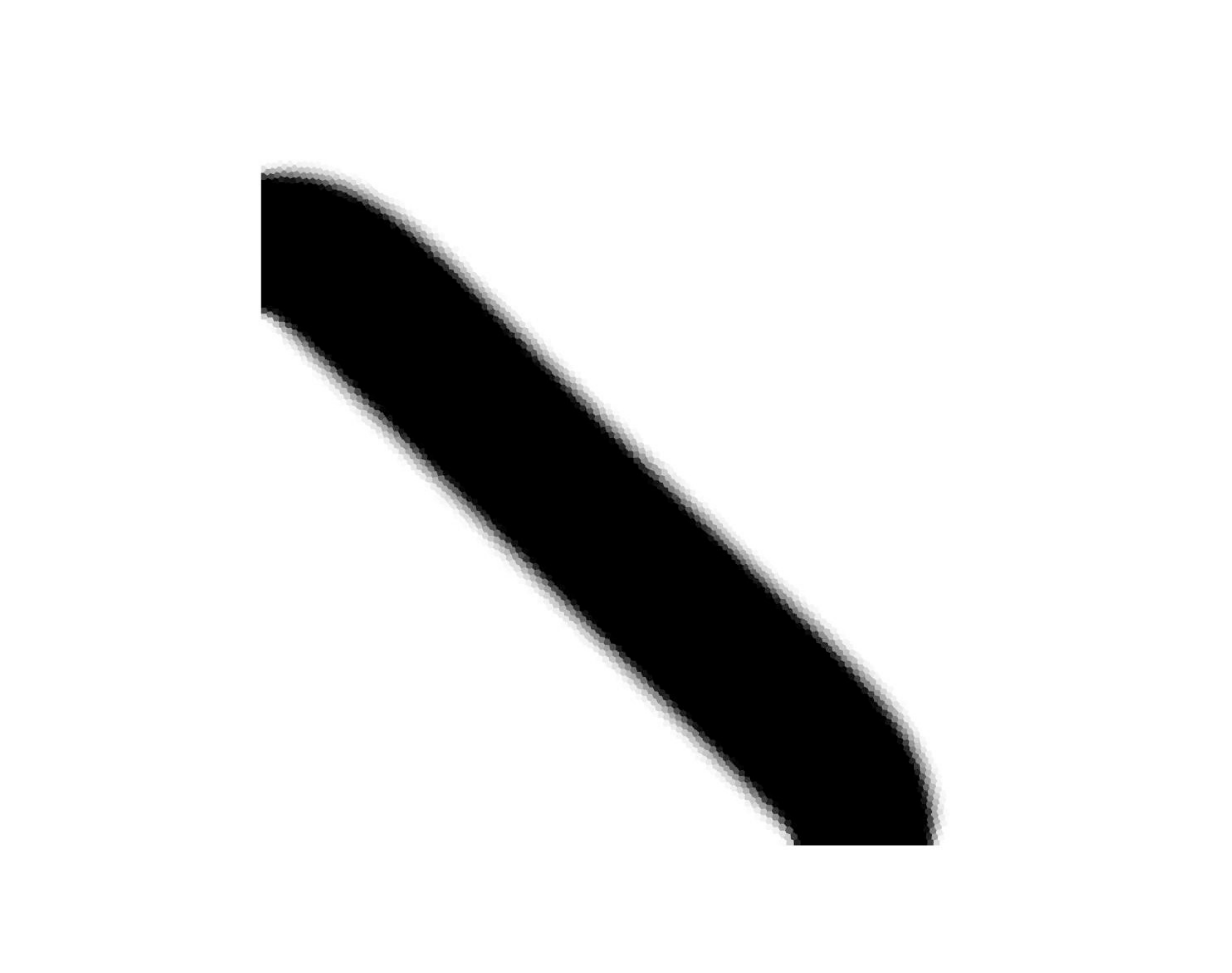}

\caption{Optimal pipe. Top: design domain (taken from  \cite{Borrvall-Petersson:03}). Bottom: optimal configuration on unstructured polygonal grid consisting of $4096$ (left) and $16384$ (right) elements.}\label{f:stokes-pipe}

\end{figure}

Here we consider the design region depicted in Figure~\ref{f:stokes-pipe}(top) where the inflow equals the outflow and set $g^*=1$. We employ an unstructured mesh of hexagons made of $4096$ and $16384$ elements, respecticely and the obtained optimal shapes of the pipe minimizing the dissipated energy are shown in  Figure~\ref{f:stokes-pipe}(bottom). The optimal result is in agreement with the configuration found in the literature (see, e.g., \cite{Borrvall-Petersson:03}). 

\subsubsection{Optimal pipe with obstacle}\label{S:pipe-obstacle}
\begin{figure}
\centering
\includegraphics[width=0.45\textwidth]{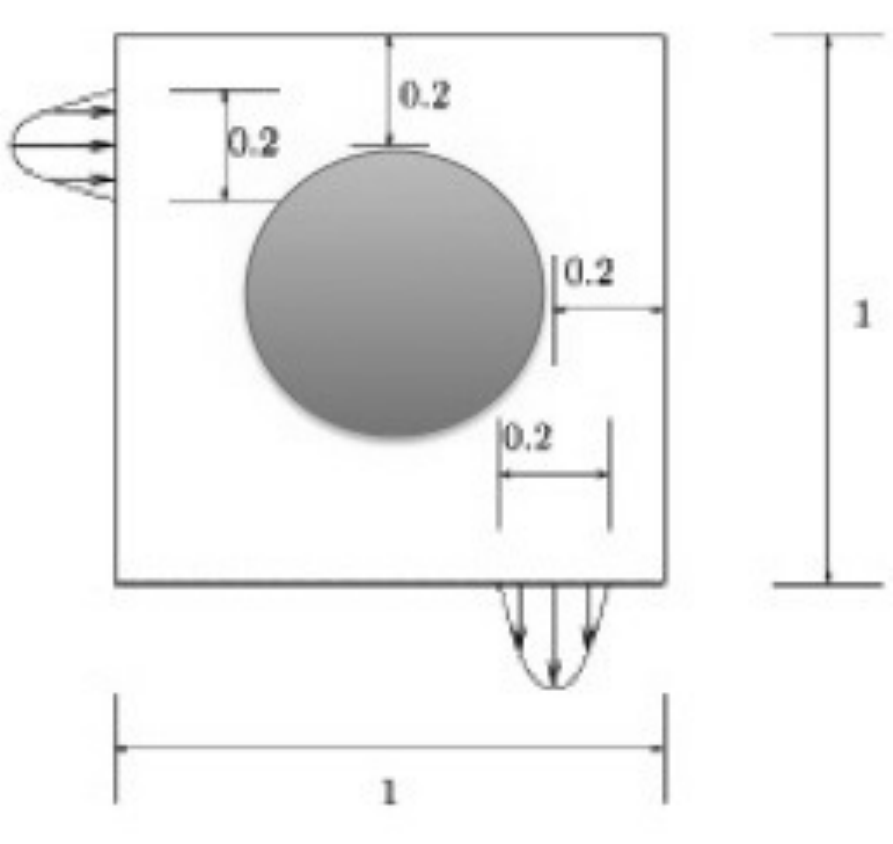}\\
\includegraphics[width=0.45\textwidth]{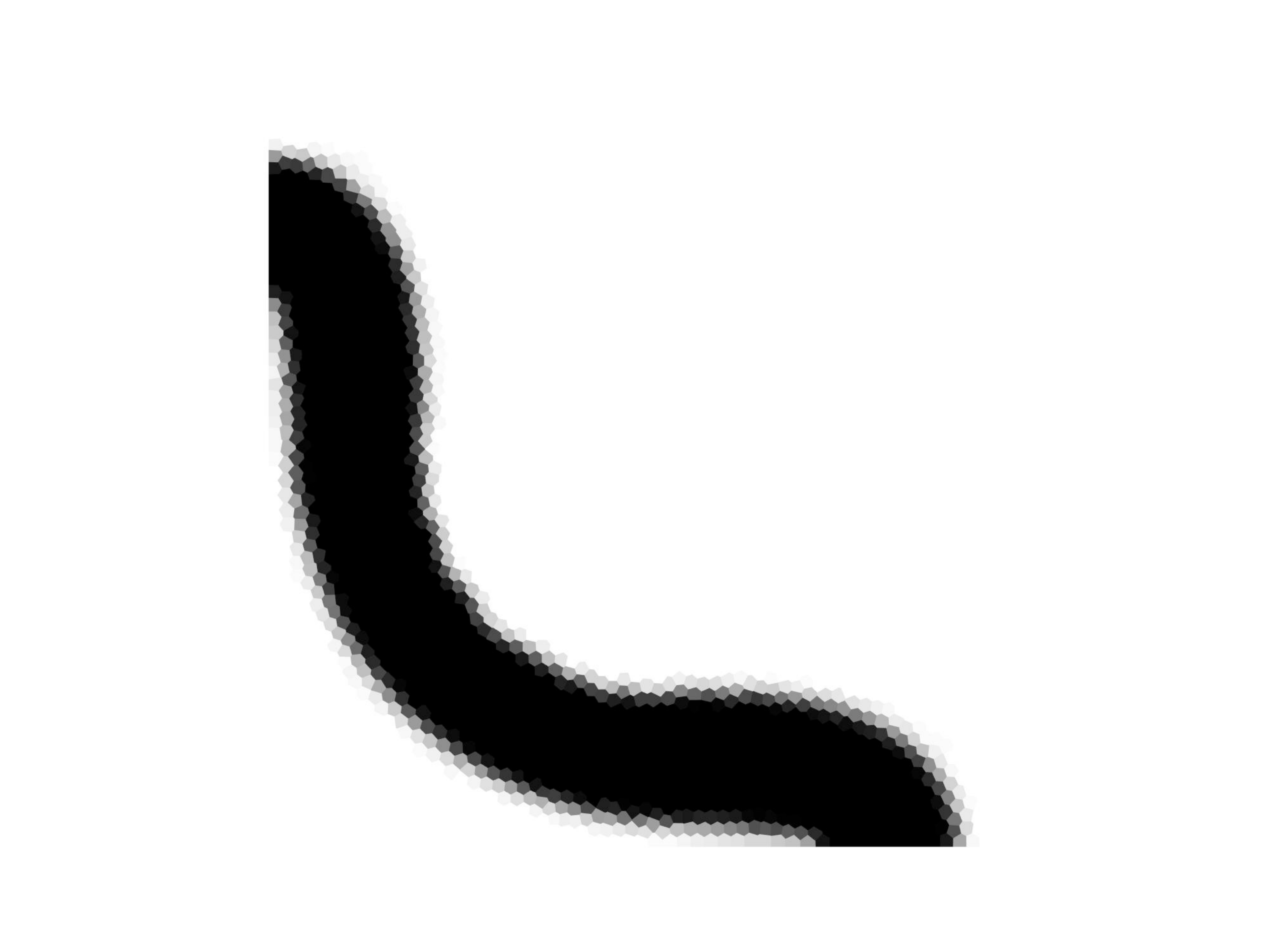}
\includegraphics[width=0.45\textwidth]{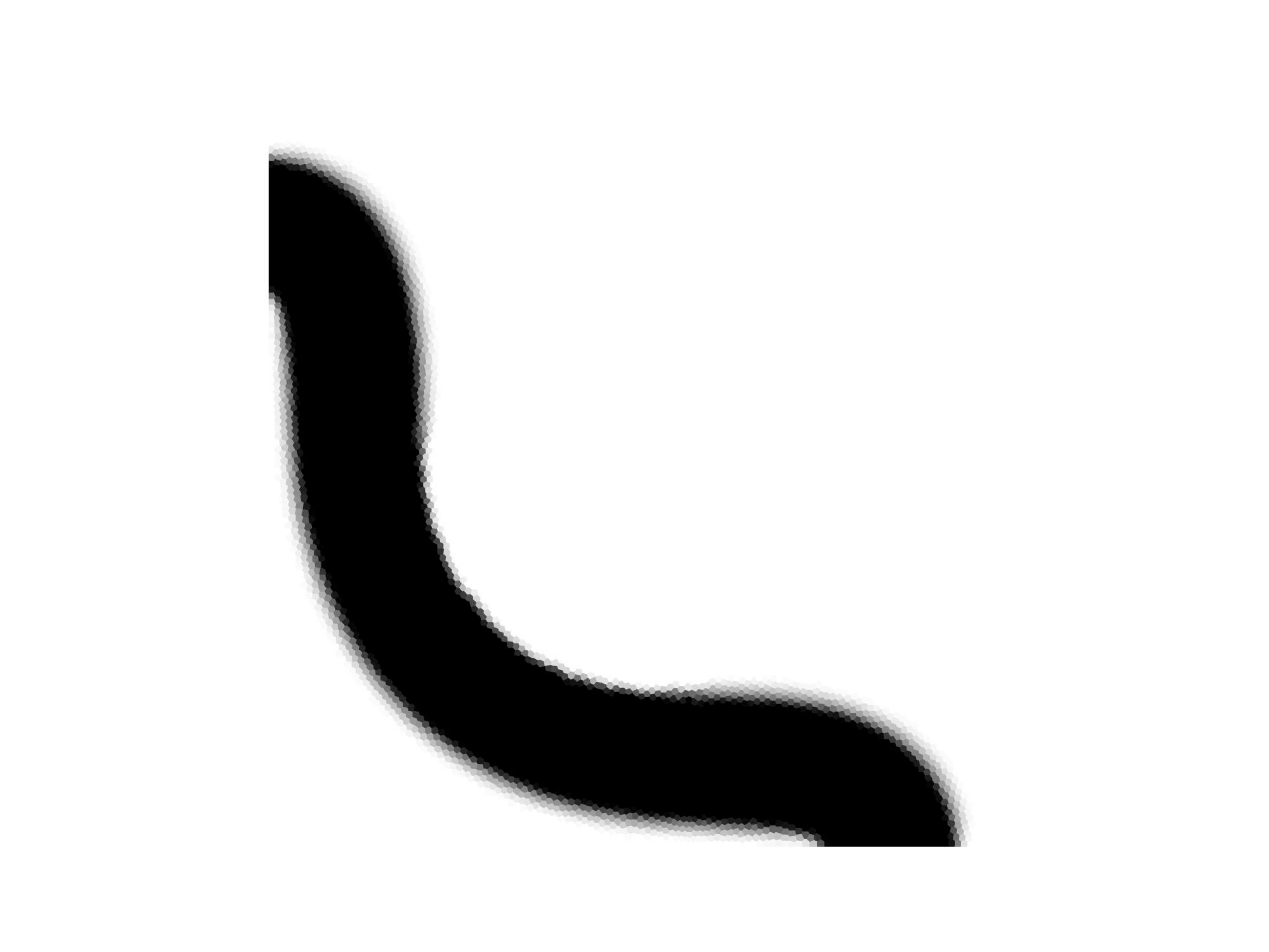}

\caption{Optimal pipe with obstacle. Top: design domain with circular obstacle. Bottom: optimal configuration on unstructured polygonal grid consisting of $4096$ (left) and $16384$ (right) elements.}\label{f:stokes-pipe-obstacle}

\end{figure}

In this section we modify the previous test case by including an obstacle represented by a circle centered at $C=(0.5,0.5)$ with radius $r=0.3$ (see Figure \ref{f:stokes-pipe-obstacle}(top)). Thus the circle is a non-design region, i.e. $\dens_h=\dens_{\textrm{min}}$, and the optimal flow has to accommodate the presence of the obstacle. 
The results of the optimization process are reported in Figure~\ref{f:stokes-pipe-obstacle}(bottom)
and clearly show the capability of the polygonal mesh to accommodate curved no-design regions.

\subsubsection{Optimal diffuser} 
\begin{figure}
\centering
\includegraphics[width=0.53\textwidth]{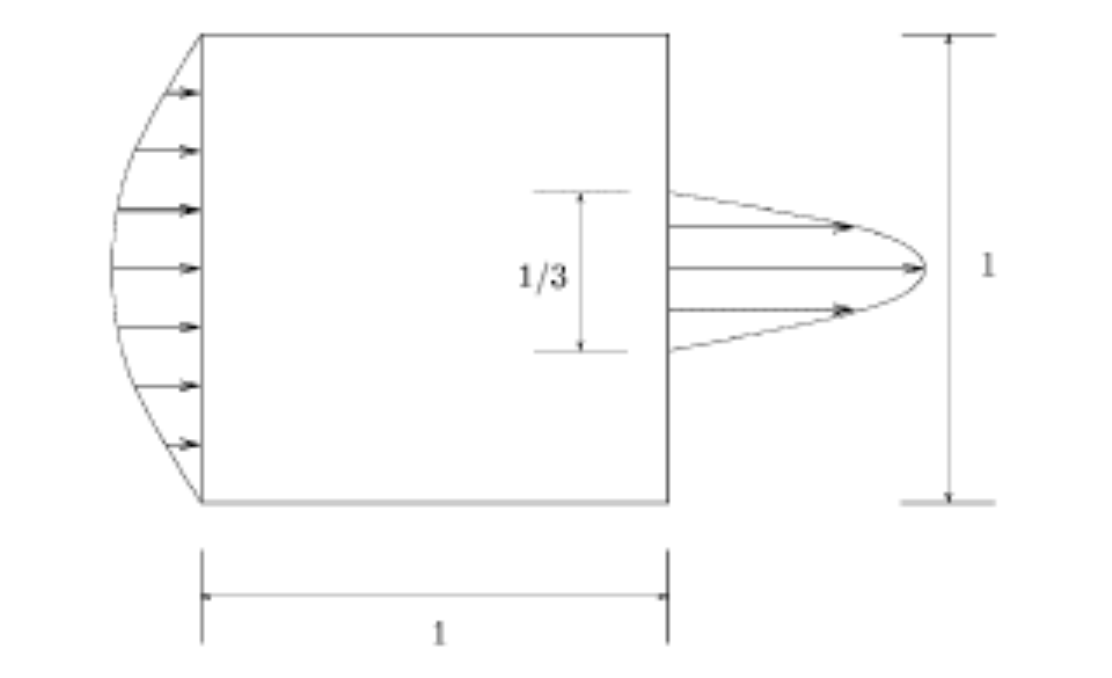}\\
\includegraphics[width=0.47\textwidth]{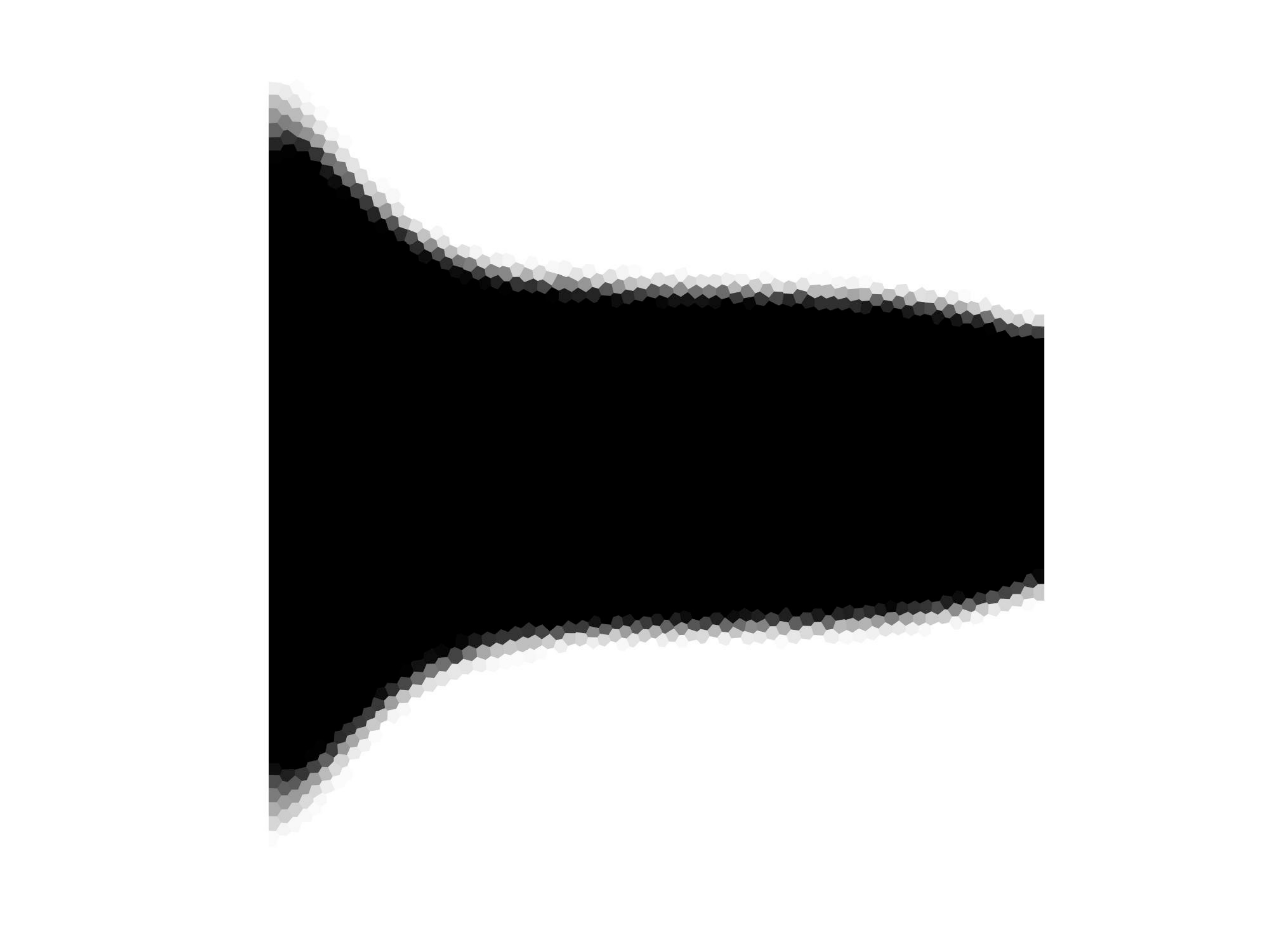}
\includegraphics[width=0.47\textwidth]{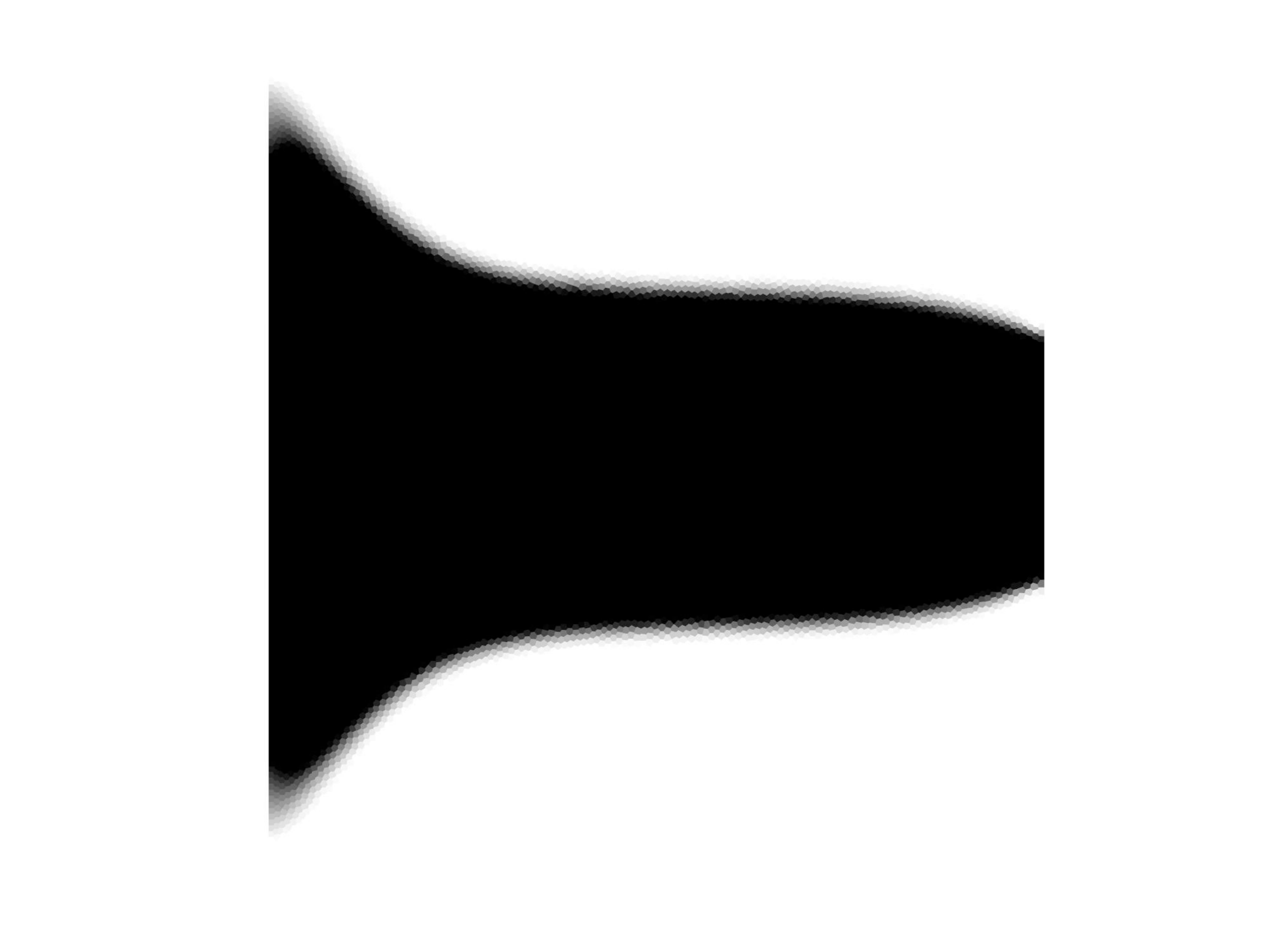}

\caption{Optimal diffuser. Top: design domain. Bottom: optimal configuration on unstructured polygonal grid consisting of $4096$ (left) and $16384$ (right) elements.}\label{f:stokes-diffuser}

\end{figure}

Here we consider the design region depicted in Figure \ref{f:stokes-diffuser} ~(top) where the inflow and the outflow have been chosen to respect the mass conservation, i.e. $g^*=1$ at the inlet and $g^*=3$ at the outlet. We employ an unstructured mesh of hexagons made of $4096$ 
and $16384$ elements, respectively, and the obtained optimal shape of the pipe minimizing the dissipated energy is shown in  Figure \ref{f:stokes-diffuser}~(bottom). Also in this case, the optimal result obtained with our VEM based approach is in agreement with the configuration found in the literature (see, e.g., \cite{Borrvall-Petersson:03}). 

\subsubsection{Optimal diffuser with obstacle} 
\begin{figure}
\centering
\includegraphics[width=0.47\textwidth]{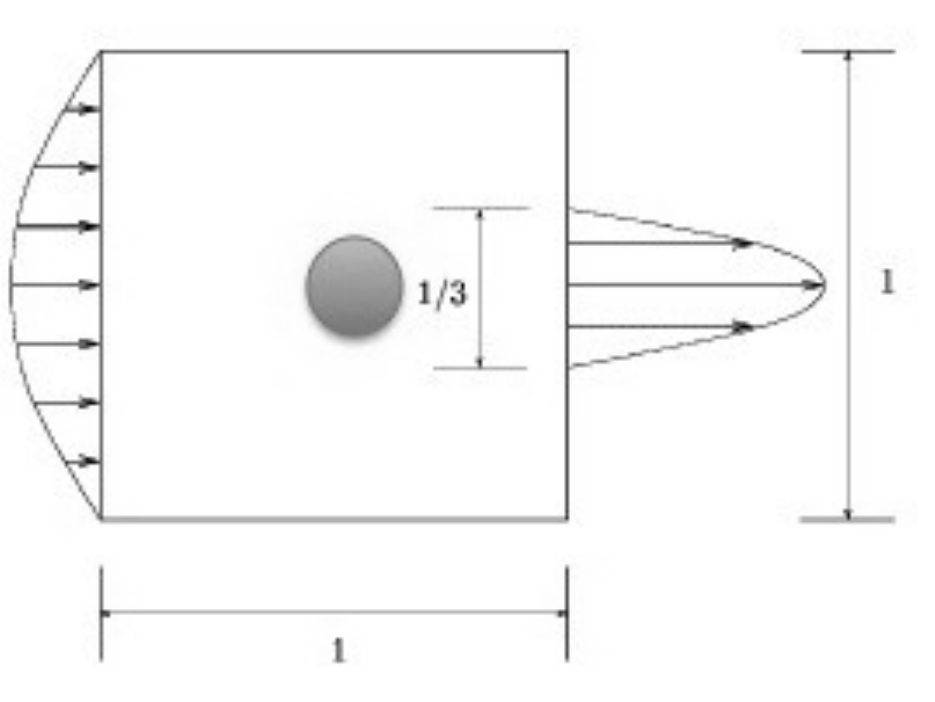} \\
\includegraphics[width=0.47\textwidth]{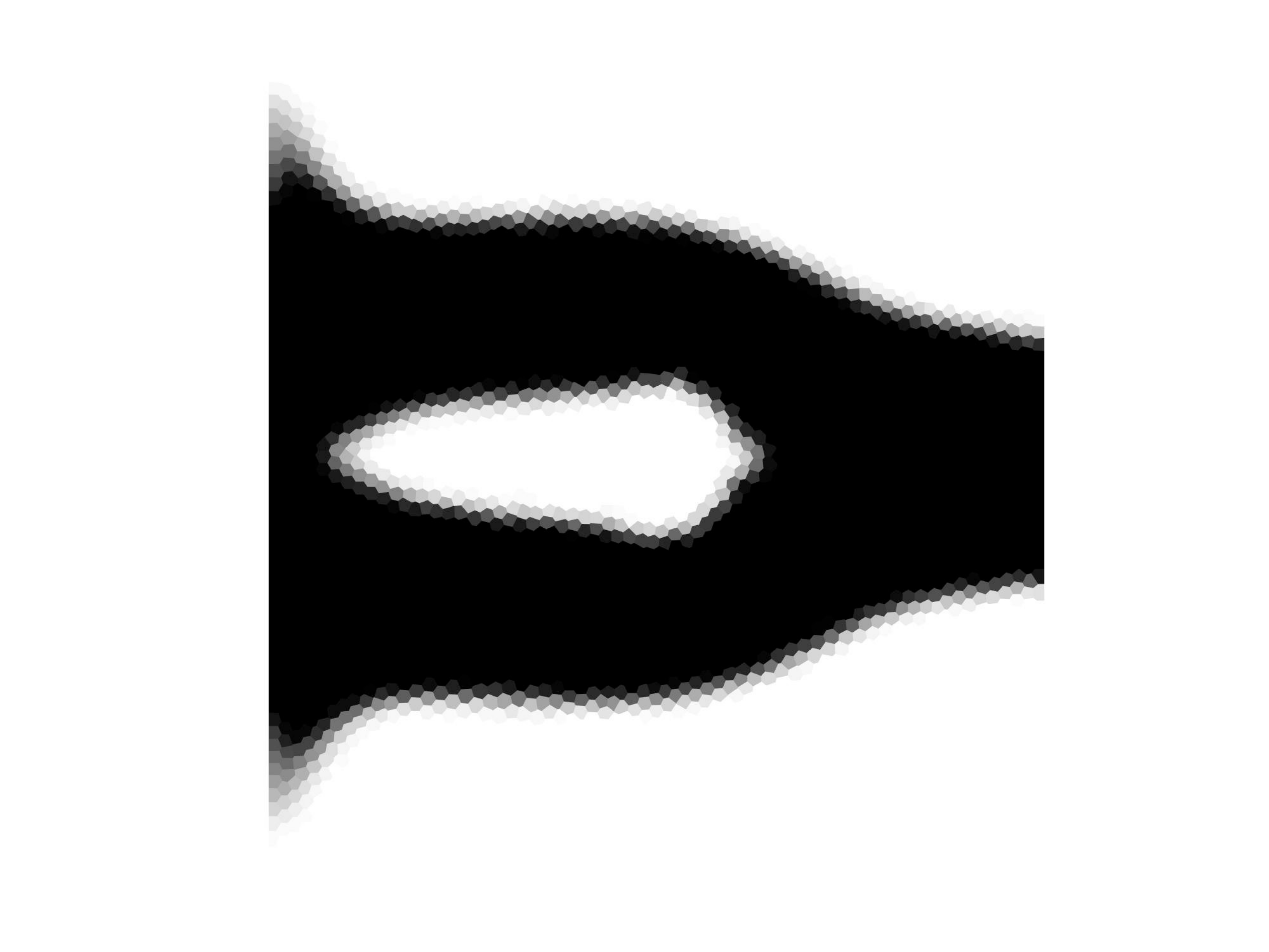}
\includegraphics[width=0.47\textwidth]{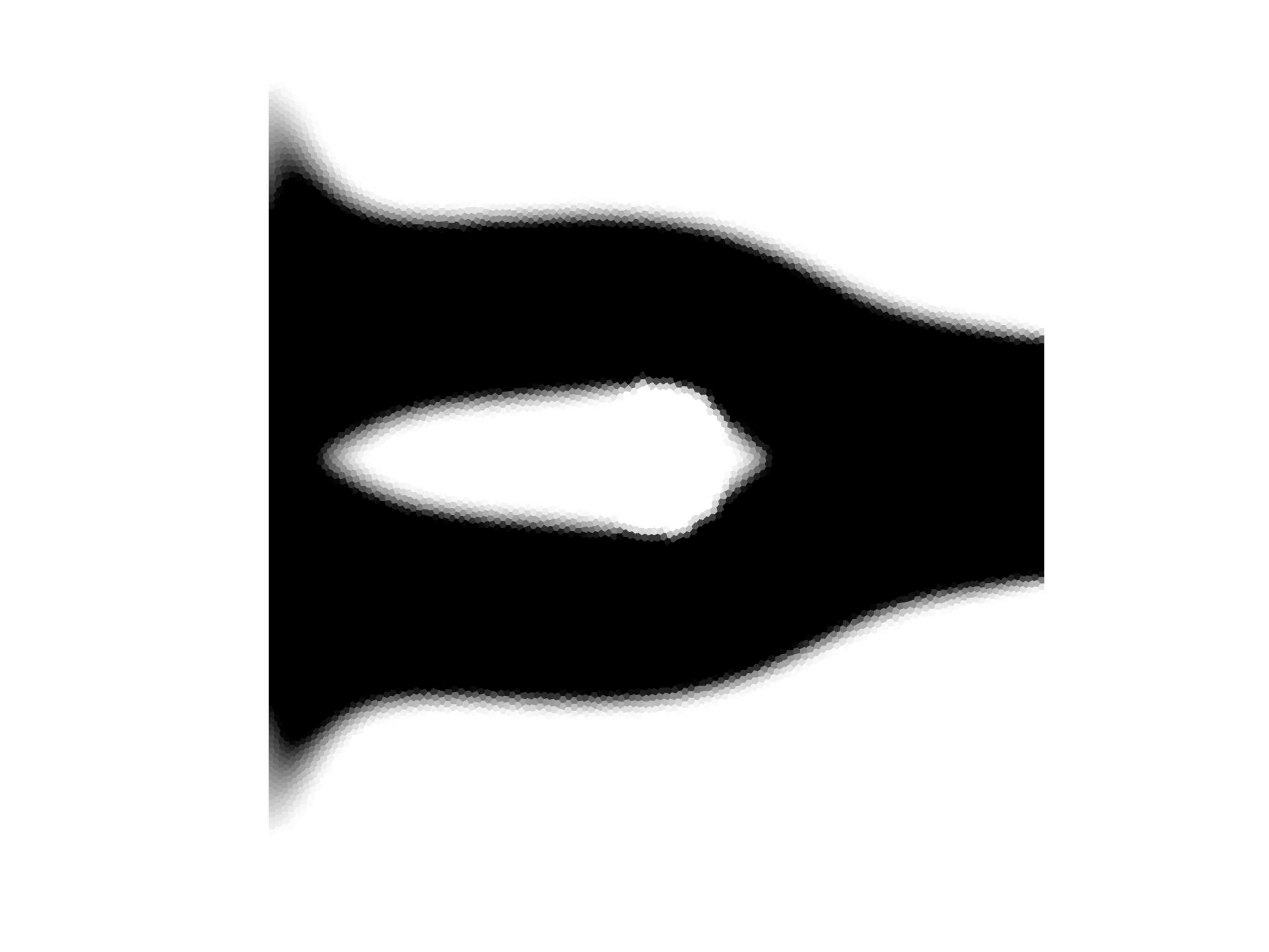}

\caption{Optimal diffuser with obstacle. Top: design domain with obstacle. Bottom: optimal configuration on polygonal grid consisting of $4096$ (left) and $16384$ (right) elements.}\label{f:stokes-diffuser-obstacle}
\end{figure}

Similarly to the optimal pipe test, we modify the previous test case by including a circular obstacle ($\dens_h=\dens_{\textrm{min}}$) centered at $C=(0.5,0.5)$ with radius $r=0.1$ (see Figure~\ref{f:stokes-diffuser-obstacle}(top)). The optimal flow taking into account the presence of the obstacle is reported in Figure~\ref{f:stokes-diffuser-obstacle}(bottom) for two polygonal grids with $4096$ and $16384$ elements, respectively. The same comments of Section \ref{S:pipe-obstacle} apply here.

\section{Conclusions}\label{S:5}
In this paper we considered the numerical solution of two paradigmatic examples of topology optimization problems on polygonal meshes employing the virtual element method. The first optimization problem is the minimum compliance governed by linear elasticity (compressible and nearly-incompressible) while the second one is related to the minimum energy dissipation of Stokes flows. From the numerical results 
presented in the previous section, we can draw the following conclusions:
\begin{itemize}
\item optimal layouts do not seem to be affected by the geometrical features of the polygonal mesh, whereas the use of standard quadrilateral grids may steer the optimization process towards sub-optimal (non-physical) configurations;
\item optimal layouts obtained with polygonal meshes seem to be robust with respect to mesh rotations, whereas rotated standard quadrilateral grids may give rise to different optimal configurations;
\item optimal configurations obtained on polygonal meshes seem to be mesh independent, i.e. they do not seem to depend on the granularity of the computational mesh;
\item optimal layouts in the case of topology optimization governed by nearly-incompressible elasticity or Stokes flow have been successfully identified, thanks to the accuracy and stability properties of the adopted virtual element approximation.  
 \end{itemize}
Hinging on the results of the present paper, it seems to be promising, in terms of reduction of the overall computational cost,  the adoption of mesh adaptivity during the optimization process (see, e.g., \cite{5,Paulettietal:2012}). This will be addressed in a future work.


\begin{thebibliography}{10}

\bibitem{Enhanced-VEM:2013}
B.~Ahmad, A.~Alsaedi, F.~Brezzi, L.~D. Marini, and A.~Russo.
\newblock Equivalent projectors for virtual element methods.
\newblock {\it Comput. Math. Appl.}, 66(3):376--391, 2013.

\bibitem{Antonietti-Beirao-Mora-Verani:2014}
P.~F. Antonietti, L.~Beir{\~a}o~da Veiga, D.~Mora, and M.~Verani.
\newblock A stream virtual element formulation of the {S}tokes problem on
  polygonal meshes.
\newblock {\it SIAM J. Numer. Anal.}, 52(1):386--404, 2014.

\bibitem{Antonietti-Beirao-Scacchi-Verani:2016}
P.~F. Antonietti, L.~Beir{\~a}o~da Veiga, S.~Scacchi, and M.~Verani.
\newblock A {$C^1$} virtual element method for the {C}ahn-{H}illiard equation
  with polygonal meshes.
\newblock {\it SIAM J. Numer. Anal.}, 54(1):34--56, 2016.

\bibitem{ncvem-Antonietti-Manzini-Verani:2016}
P. F. Antonietti, G. Manzini and M. Verani. 
\newblock The fully nonconforming Virtual Element Method for Biharmonic problems.  
\newblock {\it ArXiv e-prints}: 1611.08736, 2016.

\bibitem{Ayuso-Lipnikov-Manzini:2016}
B.~Ayuso~de Dios, K.~Lipnikov, and G.~Manzini.
\newblock The nonconforming virtual element method.
\newblock {\it ESAIM Math. Model. Numer. Anal.}, 50(3):879--904, 2016.

\bibitem{VEM-basic:2013}
L.~Beir{\~a}o~da Veiga, F.~Brezzi, A.~Cangiani, G.~Manzini, L.~D. Marini, and
  A.~Russo.
\newblock Basic principles of virtual element methods.
\newblock {\it Math. Models Methods Appl. Sci.}, 23(1):199--214, 2013.

\bibitem{Vem-elasticity:2013}
L.~Beir{\~a}o~da Veiga, F.~Brezzi, and L.~D. Marini.
\newblock Virtual elements for linear elasticity problems.
\newblock {\it SIAM J. Numer. Anal.}, 51(2):794--812, 2013.

\bibitem{Hdiv-vem:2016}
L.~Beir{\~a}o~da Veiga, F.~Brezzi, L.~D. Marini, and A.~Russo.
\newblock {$H(\text{div})$} and {$H(\bold{curl})$}-conforming virtual element
  methods.
\newblock {\it Numer. Math.}, 133(2):303--332, 2016.

\bibitem{mixed-vem:2016}
L.~Beir{\~a}o~da Veiga, F.~Brezzi, L.~D. Marini, and A.~Russo.
\newblock Mixed virtual element methods for general second order elliptic
  problems on polygonal meshes.
\newblock {\it ESAIM Math. Model. Numer. Anal.}, 50(3):727--747, 2016.

\bibitem{VEM-elliptic:2016}
L.~Beir{\~a}o~da Veiga, F.~Brezzi, L.~D. Marini, and A.~Russo.
\newblock Virtual element method for general second-order elliptic problems on
  polygonal meshes.
\newblock {\it Math. Models Methods Appl. Sci.}, 26(4):729--750, 2016.

\bibitem{Serendipity:2016}
L.~Beir{\~a}o~da Veiga, F.~Brezzi, L.D. Marini, and A.~Russo.
\newblock Serendipity nodal vem spaces.
\newblock {\it Computers and Fluids}, 2016.

\bibitem{Beirao-Chernov-Mascotto-Russo:2016}
L.~Beir{\~a}o~da Veiga, A.~Chernov, L.~Mascotto, and A.~Russo.
\newblock Basic principles of {$hp$} virtual elements on quasiuniform meshes.
\newblock {\it Math. Models Methods Appl. Sci.}, 26(8):1567--1598, 2016.

\bibitem{BeiraoErn_2016}
L.~Beir{\~a}o~da Veiga and A.~Ern.
\newblock Preface [{S}pecial issue---{P}olyhedral discretization for {PDE}].
\newblock {\it ESAIM Math. Model. Numer. Anal.}, 50(3):633--634, 2016.

\bibitem{BeiraoLipni:2010}
L.~Beir{\~a}o~da Veiga and K.~Lipnikov.
A mimetic discretization of the Stokes problem with selected edge bubbles.
\newblock {\it SIAM J. Sci. Comput.}, 32(2):875--893, 2010.

\bibitem{MFD:book}
L.~Beir{\~a}o~da Veiga, K.~Lipnikov, and G.~Manzini.
\newblock {\it The mimetic finite difference method for elliptic problems},
  volume~11 of {\it MS\&A. Modeling, Simulation and Applications}.
\newblock Springer, Cham, 2014.

\bibitem{VEM:inelastic:2015}
L.~Beir{\~a}o~da Veiga, C.~Lovadina, and D.~Mora.
\newblock A virtual element method for elastic and inelastic problems on
  polytope meshes.
\newblock {\it Comput. Methods Appl. Mech. Engrg.}, 295:327--346, 2015.

\bibitem{Beirao-Lovadina-Vacca:2015}
L.~Beir{\~a}o~da Veiga, C.~{Lovadina}, and G.~{Vacca}.
\newblock {Divergence free Virtual Elements for the Stokes problem on polygonal
  meshes}.
\newblock {\it ArXiv e-prints}: 1510.01655, October 2015.

\bibitem{Bellomo-Brezzi-Manzini:2014}
N.~Bellomo, F.~Brezzi, and G.~Manzini.
\newblock Recent techniques for pde discretizations on polyhedral meshes.
\newblock {\it Math. Models Methods Appl. Sci.}, 24:1453--1455, 2014.
\newblock (special issue).

\bibitem{4} 
M.P.~Bends{\o}e, O. Sigmund, \textit{Topology optimization theory, methods and applications}, New York, Springer, 2003.

\bibitem{SIMP1} 
M.P.~ Bends{\o}e, 
\textit{Optimal shape design as a material distribution problem}, Structural Optimization, \textbf{1} (1989), 193--202.

\bibitem{Berrone-SUPG:2016}
M.~F. Benedetto, S.~Berrone, A.~Borio, S.~Pieraccini, and S.~Scial{\`o}.
\newblock Order preserving {SUPG} stabilization for the virtual element
  formulation of advection--diffusion problems.
\newblock {\it Comput. Methods Appl. Mech. Engrg.}, 311:18--40, 2016.

\bibitem{VEM-DFN:2014}
M.~F. Benedetto, S.~Berrone, S.~Pieraccini, and S.~Scial{\`o}.
\newblock The virtual element method for discrete fracture network simulations.
\newblock {\it Comput. Methods Appl. Mech. Engrg.}, 280:135--156, 2014.

\bibitem{Borrvall-Petersson:03} 
T.~Borrvall, J.~Petersson. 
\newblock Topology optimization of fluids in Stokes flow.
\newblock{\it Internat. J. Numer. Methods Fluids},  41(1): 77--107, 2003.
		
\bibitem{Brezzi-Falk-Marini:2014}
F.~Brezzi, R.~S. Falk, and L.~Donatella~Marini.
\newblock Basic principles of mixed virtual element methods.
\newblock {\it ESAIM: Mathematical Modelling and Numerical Analysis},
  48(4):1227--1240, 2014.

\bibitem{Brezzi-Marini-plate:2013}
F.~Brezzi and L.~D. Marini.
\newblock Virtual element methods for plate bending problems.
\newblock {\it Comput. Methods Appl. Mech. Engrg.}, 253:455--462, 2013.

\bibitem{VEM-discontinuous}
F.~Brezzi and L.~D. Marini.
\newblock Virtual element and discontinuous {G}alerkin methods.
\newblock In {\it Recent developments in discontinuous {G}alerkin finite
  element methods for partial differential equations}, volume 157 of {\it IMA
  Vol. Math. Appl.}, pages 209--221. Springer, Cham, 2014.

\bibitem{Bruggi:16} 
M.~ Bruggi.
\newblock Topology optimization with mixed finite elements on regular grids
\newblock {\it Comput. Methods Appl. Mech. Engrg.}, 305: 133--153, 2016.

\bibitem{filt} 
M.~Bruggi, P.~Duysinx. 
\newblock{A stress-based approach to the optimal design of structures with unilateral behavior of material or supports}.
\newblock{\it Structural and Multidisciplinary Optimization}, 48: 311--326, 2013.
  
\bibitem{Bruggi-Venini:07} 
M.~ Bruggi, P.~Venini.
\newblock Topology optimization of incompressible media using mixed finite elements.
\newblock {\it Comput. Methods Appl. Mech. Engrg.},  196(33-34): 3151--3164, 2007. 
  
\bibitem{5} 
 M. Bruggi, M. Verani, {\it A fully adaptive topology optimization algorithm with goal--oriented error control}, Computers and Structures, \textbf{89} (2011), 1481--1493.

\bibitem{Cangiani-Gyrya-Manzini:2016}
A.~{Cangiani}, V.~{Gyrya}, and G.~{Manzini}.
\newblock {The non-conforming virtual element method for the Stokes equations}.
\newblock {\it ArXiv e-prints}: 1608.01210, August 2016.

\bibitem{Cangiani-Manzini-Sutton:2015}
A.~{Cangiani}, G.~{Manzini}, and O.~J. {Sutton}.
\newblock {Conforming and nonconforming virtual element methods for elliptic
  problems}.
\newblock {\it ArXiv e-prints}: 1507.03543, July 2015.

\bibitem{Chinosi-Marini:2016}
C.~Chinosi and L.~D. Marini.
\newblock Virtual {E}lement {M}ethod for fourth order problems:
  {$L^2$}-estimates.
\newblock {\it Comput. Math. Appl.}, 72(8):1959--1967, 2016.

\bibitem{Paulino-et-al:2014}
A.L.~Gain, C.~Talischi, G.H. Paulino.  
\newblock On the virtual element method for three-dimensional linear elasticity
 problems on arbitrary polyhedral meshes.
 \newblock{\it Comput. Methods Appl. Mech. Engrg.},  282:132--160, 2014.

\bibitem{Gain-Paulino-Duarte-Menezes:2015} 
A.L. ~Gain, G.H.~Paulino, L.S.~Duarte,  I. F. M~Menezes.
\newblock Topology optimization using polytopes.
 \newblock {\it Comput. Methods Appl. Mech. Engrg.},  293:411--430, 2015.


\bibitem{Gardini-Vacca:2016}
F.~{Gardini} and G.~{Vacca}.
\newblock {Virtual Element Method for Second Order Elliptic Eigenvalue
  Problems}.
\newblock {\it ArXiv e-prints}: 1610.03675, October 2016.

\bibitem{26a} 
S.H.~Jeong, S.H.~Park, D.H.~Choi, G.H. ~Yoon.
\newblock Topology optimization considering static failure theories for ductile and brittle materials.
\newblock {\it Comput. Struct.}, 110-111: 116--132, 2012.

\bibitem{26b} 
S.H.~Jeong, S.H.~Park, D.H.~Choi, G.H. ~Yoon. 
\newblock Toward a stress-based topology optimization procedure with indirect calculation of internal finite element information.
\newblock{\it Computers and Mathematics with Applications}, 66(6):1065--108, 2013.

\bibitem{26c} 
S.H.~Jeong, S.H.~Park, D.H.~Choi, G.H. ~Yoon. 
\newblock Fatigue and static failure considerations using a topology optimization method. 
\newblock{\it Applied Mathematical Modelling}, vol. 39(3-4): 1137--1162, 2015.

\bibitem{torto} 
C.~ Le, J.~Norato, T.E.~Bruns, C~ Ha, D.D. ~Tortorelli. 
\newblock Stress--based topology optimization for continua.
\newblock {\it Structural and Multidisciplinary Optimization}, 41:605--620, 2010.

\bibitem{Lipnikov-Manzini-Shashkov:2014}
K.~Lipnikov, G.~Manzini, and M.~Shashkov.
\newblock Mimetic finite difference method.
\newblock {\it J. Comput. Phys.}, 257(part B):1163--1227, 2014.

\bibitem{VEM-Steklov:2015}
D.~Mora, G.~Rivera, and R.~Rodr{\'{\i}}guez.
\newblock A virtual element method for the {S}teklov eigenvalue problem.
\newblock {\it Math. Models Methods Appl. Sci.}, 25(8):1421--1445, 2015.

\bibitem{Paulettietal:2012}
P.~Morin, R. H.~Nochetto, M.S.~Pauletti, M.~Verani,
\newblock Adaptive finite element method for shape optimization. 
\newblock{\it ESAIM Control Optim. Calc. Var.}, 18(4):1122--1149, 2012. 

\bibitem{Perugia-Pietra-Russo:2016}
I.~Perugia, P.~Pietra, and A.~Russo.
\newblock A plane wave virtual element method for the {H}elmholtz problem.
\newblock {\it ESAIM Math. Model. Numer. Anal.}, 50(3):783--808, 2016.

\bibitem{sigpet} O. Sigmund O, J. Petersson, {\it Numerical instabilities in topology optimization: a survey on procedures dealing with checkerboards, mesh-dependencies and local minima}, Structural Optimization, \textbf{16} (1998), 68--75.

\bibitem{7} O. Sigmund, P.M. Clausen, \textit{Topology optimization using a mixed formulation: an alternative way to solve pressure load problems}, Computer Methods in Applied Mechanics and Engineering, \textbf{196} (2017), 1874--1889.


\bibitem{mma}	
K. Svanberg, 
\textit{Method of moving asymptotes - a new method for structural optimization}, International Journal for Numerical Methods in Engineering, \textbf{24} (1897), 193--202.


\bibitem{Paulino-fluid} 
C.~Talischi, A.~Pereira, G.H.~Paulino, I.F.M ~Menezes, M.S.~Carvalho.
\newblock Polygonal finite elements for incompressible fluid flow.
 \newblock{\it Internat. J. Numer. Methods Fluids}  74(2):134--151, 2014.
		

\bibitem{PolyTop} 
C.~Talischi, G.H.~Paulino, A.~Pereira, I.F.M.~Menezes.
\newblock PolyTop: a Matlab implementation of a general topology
 optimization framework using unstructured polygonal finite element meshes.
 \newblock{\it Struct. Multidiscip. Optim.},  45(3): 329--357, 2012.
		

\bibitem{PolyMesher} 
C.~Talischi, G.H.~Paulino, A.~Pereira, I.F.M.~Menezes.
\newblock PolyMesher: a general-purpose mesh generator for polygonal elements written in Matlab.
\newblock{\it Struct. Multidiscip. Optim.},  45(3): 309--328, 2012.
		
\bibitem{Talischi-Paulino-Pereira-Menezes:2010}
C.~Talischi, G.H.~Paulino, A.~Pereira, I.F.M.~Menezes.
\newblock Polygonal finite elements for topology optimization: A unifying paradigm. 
\newblock{\it Int. J. Numer. Meth. Engng.}, 82: 671--698, 2010.

\bibitem{strand7}
{\it Theoretical background to the Straus7 finite element analysis system}, Edition 1, Strand7 Pty Ltd; 2004.

\bibitem{Vacca-Beirao:2015}
G.~Vacca and L.~Beir{\~a}o~da Veiga.
\newblock Virtual element methods for parabolic problems on polygonal meshes.
\newblock {\it Numer. Methods Partial Differential Equations},
  31(6):2110--2134, 2015.

\bibitem{vem-cinesi:2016}
J.~Zhao, S.~Chen, and B.~Zhang.
\newblock The nonconforming virtual element method for plate bending problems.
\newblock {\it Math. Models Methods Appl. Sci.}, 26(9):1671--1687, 2016.

\bibitem{6} Y. Wang, Z. Kang, Q. He, \textit{Adaptive topology optimization with independent error control for separated displacement and density fields},Computers and Structures, \textbf{135} (2014), 50--61.




\end{thebibliography}
\end{document}